\documentclass[3p]{elsarticle}

 \usepackage{graphicx}

 \usepackage{amsmath,amsfonts,amssymb}
\usepackage{amssymb}

\newcommand{\w}{\mathbf{w}}
\renewcommand{\d}{\partial}
\newcommand{\f}{\mathbf{f}}
\newcommand{\F}{\mathbf{f}}
\newcommand{\Q}{\mathbf{Q}}
\newcommand{\q}{\mathbf{q}}

\newcommand{\halb}{\frac{1}{2}}

\newcommand{\A}{{\mathbf{A}}}

\newcommand{\Id}{\mathbf{I}}

\begin{document}

\begin{frontmatter}



\title{Arbitrary-Lagrangian-Eulerian ADER-WENO Finite Volume Schemes with Time-Accurate Local Time Stepping for Hyperbolic Conservation Laws} 
%
%
\author[UNITN]{Michael Dumbser} 
\ead{michael.dumbser@unitn.it}
\address[UNITN]{Laboratory of Applied Mathematics \\ Department of Civil, Environmental and Mechanical Engineering \\ University of Trento,  Via Mesiano 77, I-38123 Trento, Italy}
\begin{abstract}
In this article a new high order accurate cell-centered Arbitrary-Lagrangian-Eulerian (ALE) Godunov-type finite volume method with time-accurate local time stepping (LTS) is 
presented. The method is by construction locally and globally conservative. The scheme is based on a one-step predictor-corrector methodology in space-time and uses three 
main building blocks: First, a high order piecewise polynomial 
WENO reconstruction, to obtain a high order data representation in space from the known cell averages of the underlying finite volume scheme. Second, a high 
order space-time Galerkin predictor step based on a weak formulation of the governing PDE on moving control volumes. Third, a high order one-step 
finite volume scheme, based directly on the integral formulation of the conservation law in space-time. 
The algorithm being entirely based on space-time control volumes naturally allows for hanging nodes also in time, hence in this framework the implementation of 
a consistent and conservative time-accurate LTS becomes very natural and simple. The method is validated on some classical shock tube problems for the Euler  
equations of compressible gasdynamics and the magnetohydrodynamics equations (MHD). The performance of the new scheme is compared with a classical high order 
ALE finite volume scheme based on global time stepping. To the knowledge of the author, this is the \textit{first} high order accurate Lagrangian finite volume 
method ever presented together with a conservative and time-accurate local time stepping feature.  
\end{abstract}
\end{frontmatter}

  
\section{Introduction}
\label{sec.intro}

After the seminal paper \cite{munz94} by Munz in 1994, many significant advances have been made in recent years concerning the construction of accurate and 
robust Godunov-type Lagrangian schemes for hydrodynamics in one and multiple space dimensions, see e.g. the cell-centered Lagrangian finite volume schemes proposed 
by Maire et al. \cite{Maire2007,Maire2010,Maire2011,MaireCyl1} and Despr\'es et al. \cite{Despres2005,Despres2009}, the staggered Lagrangian scheme presented in 
\cite{StagLag} or the very recent family of cell-centered ALE remap schemes introduced in  \cite{ShashkovCellCentered,ShashkovRemap1,ShashkovRemap2,ShashkovRemap3,ShashkovRemap4,ShashkovRemap5,MaireCyl2}. 
The first better than second order accurate cell-centered Lagrangian finite volume schemes have been proposed by Cheng et al. 
\cite{chengshu1,chengshu2,chengshu3,chengtoro} on structured grids using a nonlinear reconstruction operator of the ENO-type \cite{eno}, while the first 
better than second order accurate cell-centered Lagrangian finite volume schemes on unstructured meshes have been recently proposed by Boscheri and Dumbser in 
\cite{Lagrange2D,LagrangeNC}. For high order unstructured Lagrangian finite-element schemes the reader is referred to \cite{scovazzi1,scovazzi2} and references 
therein. 

A common shortcoming in all the above mentioned methods is the use of a global time stepping scheme, i.e. the smallest cell usually dictates the timestep of 
all control volumes in the entire computational domain. Since strong mesh deformation and cell distorsion at shear waves and massive cell clustering at shock waves 
are common features of all Lagrangian methods, and which naturally also mimick the flow physics on the discrete level, all these methods suffer from very small time 
steps compared to classical Eulerian finite volume schemes on a fixed mesh, where the mesh quality can be controlled by the user \textit{a priori}. 

In recent years, several successful attempts have been reported in literature to construct high order Eulerian schemes on \textit{fixed meshes} that allow for 
time-accurate \textit{local} time stepping (LTS), where each element can run at its own optimal time step, given by a \textit{local} CFL stability condition. 
Most of these schemes were of the discontinuous Galerkin finite element type \cite{FlahertyLTS,dumbserkaeser06d,TaubeMaxwell,stedg1,stedg2,KrivodonovaLTS}, but also high order 
accurate finite volume schemes with time accurate LTS can be found, mostly in the context of adaptive mesh refinement (AMR) methods and block-clustered local 
time-stepping, see \cite{DomainDecomp,CastroLTS,MuletAMR1,MuletAMR2,Burger2012,AMRCL,AMRNC}. High order Runge-Kutta time integrators with time-accurate 
local-time stepping have been recently proposed in \cite{GroteLTS1,GroteLTS2}. 

In this article, a first attempt is made to introduce time-accurate local time stepping also into high order cell-centered \textit{Lagrangian} finite volume schemes, in order to 
reduce the computational effort in the presence of small and highly deformed cells caused by the Lagrangian framework. The paper is organized as follows: in Section 
\ref{sec.algo} the time-accurate local time stepping algorithm is described and all important implementation details are given; in Section \ref{sec.tests} some 
numerical test problems are presented in order to validate the accuracy, efficiency and robustness of the proposed approach; detailed comparisons with global time 
stepping (GTS) are made; the paper is closed by concluding remarks and an outlook to future extensions in Section \ref{sec.concl}. 

\section{Algorithm Description}
\label{sec.algo}

In this paper, one-dimensional hyperbolic systems of conservation laws of the form 
\begin{equation}
\label{eqn.pde}
 \frac{\partial }{\partial t} \Q + \frac{\partial }{\partial x} \f(\Q) = 0, \qquad x \in \Omega(t) \subset \mathbb{R}, \quad t \in \mathbb{R}_0^+, 
\end{equation} 
are considered, where $\Q \in \Omega_Q \subset \mathbb{R}^\nu$ is the vector of conserved variables and $\f=\f(\Q)$ is the nonlinear flux vector.  
Here, $\Omega(t)$ denotes the time-dependent computational domain and $\Omega_Q$ is the set of admissible states, the so-called phase-space 
or state-space. 
The computational domain $\Omega(t)$ is discretized by a set of moving mesh points $x_{i+\halb}(t)$, each of which can move  
with a generic local mesh velocity $V_{i+\halb}(t)$, so that the trajectory of the point satisfies the ODE 
\begin{equation} 
\label{eqn.mesh.ode} 
  \frac{d}{dt} x_{i+\halb} = V_{i+\halb}(t). 
\end{equation} 
The Arbitrary-Eulerian-Lagrangian flux vector and its Jacobian are introduced as  
\begin{equation}
\F^V(\Q,V) = \F(\Q) - V \Q, \qquad \textnormal{ and } \qquad \A^V(\Q,V) = \frac{\partial \F^V}{\partial \Q}, 
\end{equation}
where $V$ is the local mesh velocity. 

The family of high order one-step cell-centered Lagrangian finite volume schemes described in this section proceeds for each element and (local) time step with the 
following three subteps: 
\begin{enumerate}
\item piecewise polynomial data reconstruction from known cell averages, $\Q_i^n \to \w_h(x,t^n)$,  
\item element-local data evolution in time, $\w_h(x,t^n) \to \q_h(x,t)$ and 
\item one-step element update $\Q_i^n \to \Q_i^{n+1}$ and geometry update 
$x_{i  \pm \halb}^n \to x_{i \pm \halb}^{n+1}$. 
\end{enumerate} 
\begin{flushright}

\end{flushright}
While this is the natural order of the 
three steps in the computer program and in the numerical scheme, in the following subsections the three steps of the scheme are described in the opposite order 
for the sake of clarity, i.e. first the finite volume scheme is presented, then the local data evolution is described and finally the data reconstruction 
step is outlined. 

\subsection{One-Step Cell-Centered Lagrangian Finite Volume Scheme with LTS} 
\label{sec.onestep} 

At time $t$ the spatial control volumes are defined as $T_i=T_i(t) = [x_{i-\halb}(t); x_{i+\halb}(t)]$. The aim is first to construct a 
one-step finite volume scheme in space and time, because the use of a one-step time discretization is one of the key ingredients to achieve a simple and efficient 
local time stepping algorithm\footnote{In contrast to the usual Runge-Kutta time stepping, which requires several sub-stages.}. 
Following \cite{Lagrange1D}, the conservation law \eqref{eqn.pde} is first integrated in space and time over a generic space-time control volume 
$\mathcal{C}_i = T_i(t) \times [t_i^{n}; t_i^{n+1}]$, 
\begin{equation}
\label{eqn.pde.int} 
  \Delta x_i^{n+1} \Q_i^{n+1} = \Delta x_i^{n} \Q_i^{n} - 
	\left( \int \limits_{t_{i}^n}^{t_{i}^{n+1}}  \F^V\left(\Q(x_{i+\halb}(t),t), V_{i+\halb}(t) \right)   dt - 
	       \int \limits_{t_{i}^n}^{t_{i}^{n+1}}  \F^V\left(\Q(x_{i-\halb}(t),t), V_{i-\halb}(t) \right)   dt \right). 			     
\end{equation} 
The integral formulation \eqref{eqn.pde.int} is by construction locally and globally conservative. 
If one wants to achieve time-accurate local time stepping (LTS) then a cell $T_i$ has to satisfy the so-called \textit{update criterion} 
or \textit{evolve condition} \cite{dumbserkaeser06d,stedg1}   
\begin{equation}
 \label{eqn.update.criterion} 
    t_i^n + \Delta t_i^n \leq \min_{j \in \mathcal{N}_i} \left( t_j^n + \Delta t_j^n \right),  \quad  \textnormal{ or equivalently  } \quad 
		t_i^{n+1} \leq \min_{j \in \mathcal{N}_i} \left( t_j^{n+1} \right),
\end{equation} 
where $\mathcal{N}_i=\{i-1,i+1\}$ are the side neighbors of $T_i$. The current time, the future time and the time step 
in element $T_i$ are denoted by $t_i^n$, $t_i^{n+1}$ and $\Delta t_i^n$, respectively. 
Condition \eqref{eqn.update.criterion} means that a cell can be only updated if its future time is less or equal than all the future times of the neighbor elements. 
Since a time-accurate local time-stepping (LTS) algorithm produces hanging nodes in time at the edges of an element, 
the flux integrals appearing in \eqref{eqn.pde.int} are conveniently computed at the aid of a \textit{memory variable} that properly takes into account all 
fluxes through the element interfaces $x_{i\pm\halb}$ in the \textit{past}. The finite volume scheme \eqref{eqn.pde.int} with LTS then reads   
\begin{equation}
\label{eqn.fv} 
  \Delta x_i^{n+1} \Q_i^{n+1} = \Delta x_i^{n} \Q_i^{n} - \left( \Delta t_{i+\halb}^n f^V_{i+\halb} - \Delta t_{i-\halb}^n f^V_{i-\halb} \right) + \Q_i^M.      
\end{equation} 
In the relations above $\Delta x_i(t) = x_{i+\halb}(t) - x_{i-\halb}(t)$ denotes the mesh spacing at a general time $t$. While $\Delta t_i^n = t_i^{n+1}-t_i^n$ 
is the size of the current element-local time step the $\Delta t_{i \pm \halb}^n$ denote the lengths of the time-intervals on the edges. 
These edge time-intervals (or flux time-intervals) are defined as 
\begin{equation}
\label{eqn.flux.timeint}
\Delta t_{i \pm \halb}^n = t_{i \pm \halb}^{n+1} - t_{i \pm \halb}^{n},  \quad \textnormal{ with }  
\left[t_{i \pm \halb}^{n}; t_{i \pm \halb}^{n+1}\right] = \left[\max\left( t_i^n, t_{i \pm 1}^n \right); \min\left( t_i^{n+1}, t_{i \pm 1}^{n+1} \right)  \right]. 
\end{equation} 
For $t_i^n$ the time step number $n$ is an element-local index and for $t_{i+\halb}^n$ it is an edge-local number,  
but to ease notation we simply always write $n$, intending $n=n(i)$ inside elements and $n=n(i+\halb)$ at edges. 
From \eqref{eqn.update.criterion} and \eqref{eqn.flux.timeint} it follows that $t_{i \pm \halb}^{n+1} = t_i^{n+1}$.  
Furthermore, $\Delta x_i^n = \Delta x_i(t_i^n)$ and $\Q_i^{n} = \Q_i(t_i^n)$. 
The cell averages are defined as usual as 
\begin{equation}
\label{eqn.avg} 
  \Q_i(t) = \frac{1}{\Delta x_i(t)} \int \limits_{x_{i-\halb}(t)}^{x_{i+\halb}(t)} \Q(x,t) dx,    
\end{equation}  
while the time-averaged ALE interface flux across the element boundaries is defined as  
\begin{equation}
\label{eqn.fluxdef} 
  \F^V_{i+\halb} = \frac{1}{\Delta t_{i+\halb}^n} \int \limits_{t_{i+\halb}^n}^{t_{i+\halb}^{n+1}}  \F^V\left(\Q(x_{i+\halb}(t),t), V_{i+\halb}(t) \right)  \ dt.    
\end{equation}
In \eqref{eqn.fv} $\Q_i^M$ is the \textit{memory variable} \cite{dumbserkaeser06d} that takes into account all fluxes through the element interfaces $x_{i\pm \halb}$ 
in the past between time $t_i^n$ and the times $t_{i\pm\halb}^n$, i.e. 
\begin{equation}
\label{eqn.memvar.def} 
  \Q_i^M = - 
	\left( \int \limits_{t_{i}^n}^{t_{i+\halb}^n}  \F^V\left(\Q(x_{i+\halb}(t),t), V_{i+\halb}(t) \right)   dt - 
	       \int \limits_{t_{i}^n}^{t_{i-\halb}^n}  \F^V\left(\Q(x_{i-\halb}(t),t), V_{i-\halb}(t) \right)   dt \right). 
\end{equation} 
Using the definition of the memory variable \eqref{eqn.memvar.def} and the LTS finite volume scheme \eqref{eqn.fv} with \eqref{eqn.fluxdef} one obtains again the 
original integral form of the conservation law \eqref{eqn.pde.int}.  
In practice, however, the memory variable is \textit{not} computed directly by formula \eqref{eqn.memvar.def}, but according to the following strategy, see also  
\cite{dumbserkaeser06d}:  
After a local time-step has been carried out by element $T_i$ according to \eqref{eqn.fv}, then its memory variable is reset to zero and the fluxes 
computed through the element interfaces are immediately \textit{accumulated} into the memory variables of the \textit{neighbor cells} (to assure conservation), i.e. 
\begin{equation}
\label{eqn.memvar}
 \Q_i^M := 0, \qquad \Q_{i \pm 1}^M := \Q_{i \pm 1}^M \pm \Delta t_{i \pm \halb}^n f^V_{i \pm \halb}. 
\end{equation}  
In other words, the flux contributions to the memory variable of a cell $i$ are always computed by the \textit{neighbor} elements and only the 
reset step is done by the element itself. It is easy to see that algorithm \eqref{eqn.memvar} is equivalent with definition \eqref{eqn.memvar.def}, since
the time integrals in \eqref{eqn.memvar.def} are additive. 

While Eqn. \eqref{eqn.fv} with \eqref{eqn.avg}, \eqref{eqn.fluxdef} and \eqref{eqn.memvar.def} is an \textit{exact} integral relation, a numerical scheme is 
obtained by  using a \textit{numerical flux} $\F^V_h(\q_h^-,\q_h^+,V)$ instead of \eqref{eqn.fluxdef}, where the flux becomes a function of \textit{two} state 
vectors,  namely the states $\q_h^- = \q_h(x_{i+\halb}^-(t),t)$ and $\q_h^+ = \q_h(x_{i+\halb}^+(t),t)$ on the left and on the right of the interface, respectively,  
\begin{equation} 
\label{eqn.numflux} 
\F^V_{i+\halb} = \frac{1}{\Delta t_{i+\halb}^n} \int \limits_{t_{i+\halb}^n}^{t_{i+\halb}^{n+1}} \F^V_h \left(\q_h(x_{i+\halb}^-(t),t),\q_h(x_{i+\halb}^+(t),t), V_{i+\halb}(t) \right) dt.  
\end{equation} 
The procedure for the computation of $\q_h(x,t)$ will be described in the next section.   
In this article, two different numerical fluxes are used: either a simple Rusanov-type flux \cite{Rusanov:1961a}, or an Osher-type flux, 
as introduced in \cite{OsherUniversal,OsherNC}. The Rusanov-type flux reads 
\begin{equation}
\label{eqn.rusanov} 
  \F^V_h(\q_h^-,\q_h^+,V_{i+\halb}) = \frac{1}{2} \left( \F^V(\q_h^-,V_{i+\halb}) + \F^V(\q_h^+,V_{i+\halb}) \right) - \frac{1}{2} s_{\max} \left( \q_h^+ - \q_h^- \right),   
\end{equation} 
where $s_{\max} = \max( \max( |\lambda(\A^V(\q_h^-,V_{i+\halb}))|), \max( |\lambda(\A^V(\q_h^+,V_{i+\halb}))|) $ is the maximum signal speed.  
The Osher-type flux according to \cite{OsherUniversal,Lagrange1D} reads 
\begin{equation}
\label{eqn.osher} 
  \F^V_h(\q_h^-,\q_h^+,V_{i+\halb}) = \frac{1}{2} \left( \F^V(\q_h^-,V_{i+\halb}) + \F^V(\q_h^+,V_{i+\halb}) \right) - \frac{1}{2} \left( \int \limits_0^1 \left|\A^V( \mathbf{\Psi}(s), V_{i+\halb})\right| ds \right) \left( \q_h^+ - \q_h^- \right),   
\end{equation} 
where 
\begin{equation}
  \mathbf{\Psi}(s) = \mathbf{\Psi}(\q_h^-,\q_h^+,s) = \q_h^- + s \left(\q_h^+ - \q_h^- \right) 
\end{equation} 
is a straight-line segment path connecting the two states $\q_h^-$ and $\q_h^+$ in phase-space. The integral appearing in Eqn. \eqref{eqn.osher} is computed 
\textit{numerically} using Gauss-Legendre quadrature formulae of appropriate order, see \cite{OsherUniversal,OsherNC}.  
In \eqref{eqn.osher}, the usual definition for the absolute value of a matrix $\A$ applies: 
\begin{equation}
 | \A | = \mathbf{R} |\mathbf{\Lambda}| \mathbf{R}^{-1}, 
\end{equation} 
with $\mathbf{R}$ the matrix of right-eigenvectors and $\mathbf{R}^{-1}$ its inverse, while $|\mathbf{\Lambda}|=\textnormal{diag}(|\lambda_1|,|\lambda_2|,...,|\lambda_\nu|)$ is the diagonal matrix of the absolute values of the eigenvalues of $\A$. 

For the mesh velocity, needed in \eqref{eqn.meshmove} and in the fluxes \eqref{eqn.rusanov} and \eqref{eqn.osher} the Roe averaged velocity 
for Lagrangian  gasdynamics is used, see \cite{munz94},       
\begin{equation}
\label{eqn.meshvelocity}  
V_{i+\halb} = \frac{1}{2} \left( V(\q_h^-) + V(\q_h^+) \right), 
\end{equation}  
where $V = V(\Q)$ is the local fluid velocity computed from the vector of conserved variables. Note that in Lagrangian gas dynamics, the Roe average for the velocity 
is simply given by the arithmetic average of the velocities, see \cite{munz94} for details. 

Finally, using \eqref{eqn.mesh.ode} and \eqref{eqn.meshvelocity}, the new positions of the mesh points $x_{i \pm \halb}$ at times $t_{i\pm\halb}^{n+1}$ read  
\begin{equation}
\label{eqn.meshmove} 
x_{i\pm\halb}^{n+1} = x_{i\pm\halb}^n + \frac{1}{2} \int \limits_{t_{i\pm\halb}^{n}}^{t_{i\pm\halb}^{n+1}} \left( V(\q_h^-) + V(\q_h^+) \right) dt, 
\end{equation} 
which is the integral form of the ODE \eqref{eqn.mesh.ode}. 
For Lagrangian LTS schemes, it is very important to distinguish between the \textit{local element times} $t_i^{n}$ within the space-time element $\mathcal{C}_i$ and the 
\textit{local node times} $t_{i \pm \halb}^{n}$ at the nodes of the element. 
Note further that if cell number $i$ is updated, only the meshpoints $x_{i\pm\halb}$ move, the rest of the mesh remains fixed. 
In other words: if an element $T_i$ is updated in time according to \eqref{eqn.fv} it pushes the two nodes that compose element 
$T_i$ according to Eqn. \eqref{eqn.meshmove}. 

The integral conservation equation \eqref{eqn.fv} only produces one evolution equation for the cell averages $\Q_i(t)$, but since the interface flux $\F^V_{i+\halb}$ needs values at 
the element boundaries, a spatial reconstruction operator is necessary to produce appropriate interface values from the given cell averages. While a first order Godunov finite volume 
scheme simply uses extrapolation of piecewise constant data to cell boundary $x_{i+\halb}$, 
\begin{equation}
\q_h^- = \q_h(x_{i+\halb}^-,t) = \Q_i^n, \qquad \textnormal{and} \qquad \q_h^+ = \q_h(x_{i+\halb}^+,t) = \Q_{i+1}^n, 
\end{equation} 
higher order in space can be achieved in the finite volume context by using an appropriate reconstruction or recovery operator. 
In this paper, a particular form of WENO reconstruction \cite{shu_efficient_weno} is used, see \cite{Lagrange1D} and Section \ref{sec.weno}. 

\subsection{High Order Time-Evolution} 
\label{sec.time} 

The high order element-local data evolution stage is another key ingredient for the design of efficient and simple time-accurate local time stepping schemes, 
see \cite{dumbserkaeser06d,TaubeMaxwell,stedg1,stedg2}. 
Instead of the Cauchy-Kovalewski procedure, which is based on cumbersome Taylor series and repeated differentiation of the governing PDE and which has been used in the original ENO 
method of Harten et al. \cite{eno} and in the ADER schemes of Titarev and Toro \cite{toro3,titarevtoro} as well as in the Lagrangian finite volume schemes presented in \cite{chengshu2,chengtoro}, 
here a \textit{weak integral formulation} of the governing PDE in space-time is used. This concept has been introduced in \cite{DumbserEnauxToro,DumbserZanotti,HidalgoDumbser}, 
is capable of dealing with \textit{stiff} algebraic source terms and has already been successfully applied to Lagrangian finite volume schemes with global time stepping (GTS) in 
\cite{Lagrange1D,Lagrange2D,LagrangeNC}. In this section we summarize the description already given in \cite{Lagrange1D}, but for the sake of clarity the main steps are explained again 
to make this article self-contained. 

\begin{figure}[!t]
\begin{center}
\begin{tabular}{lr}
\includegraphics[width=0.45\textwidth]{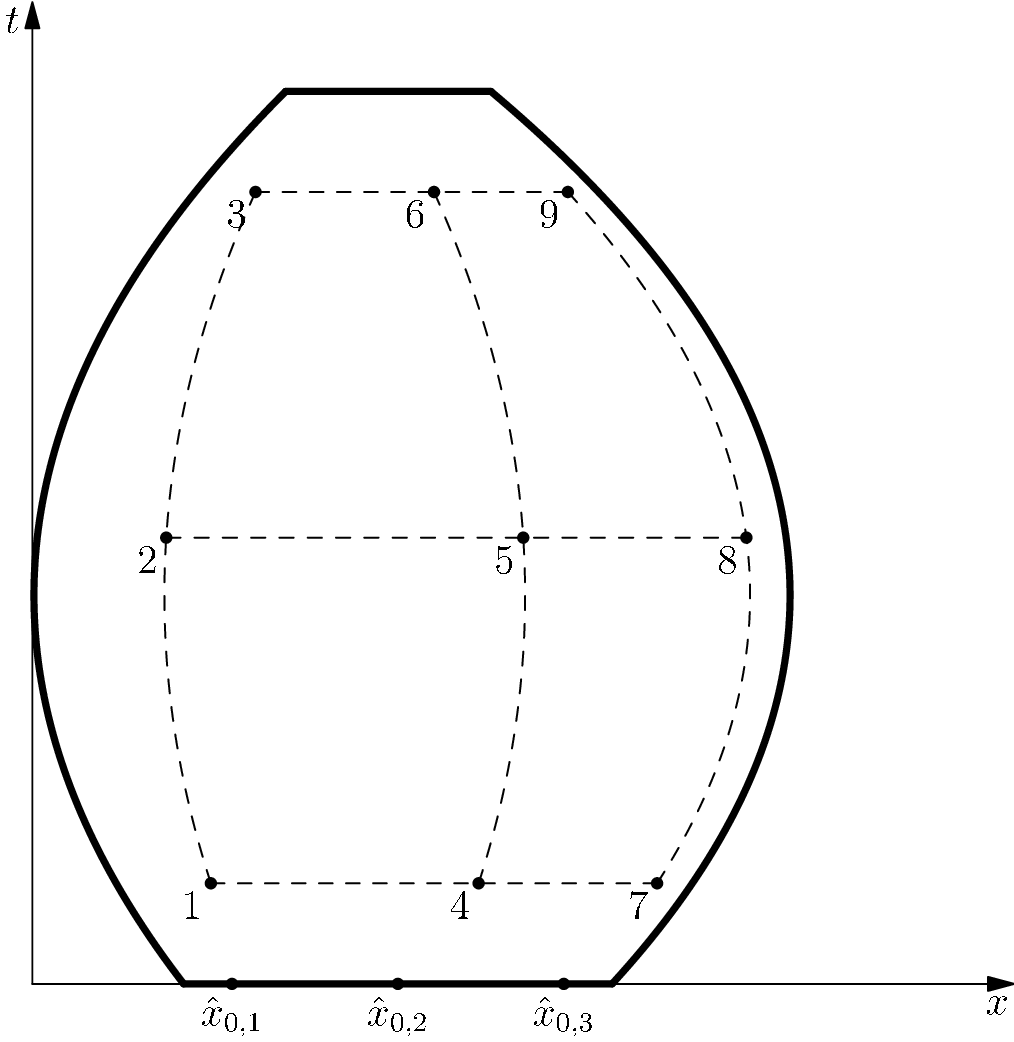} &
\includegraphics[width=0.45\textwidth]{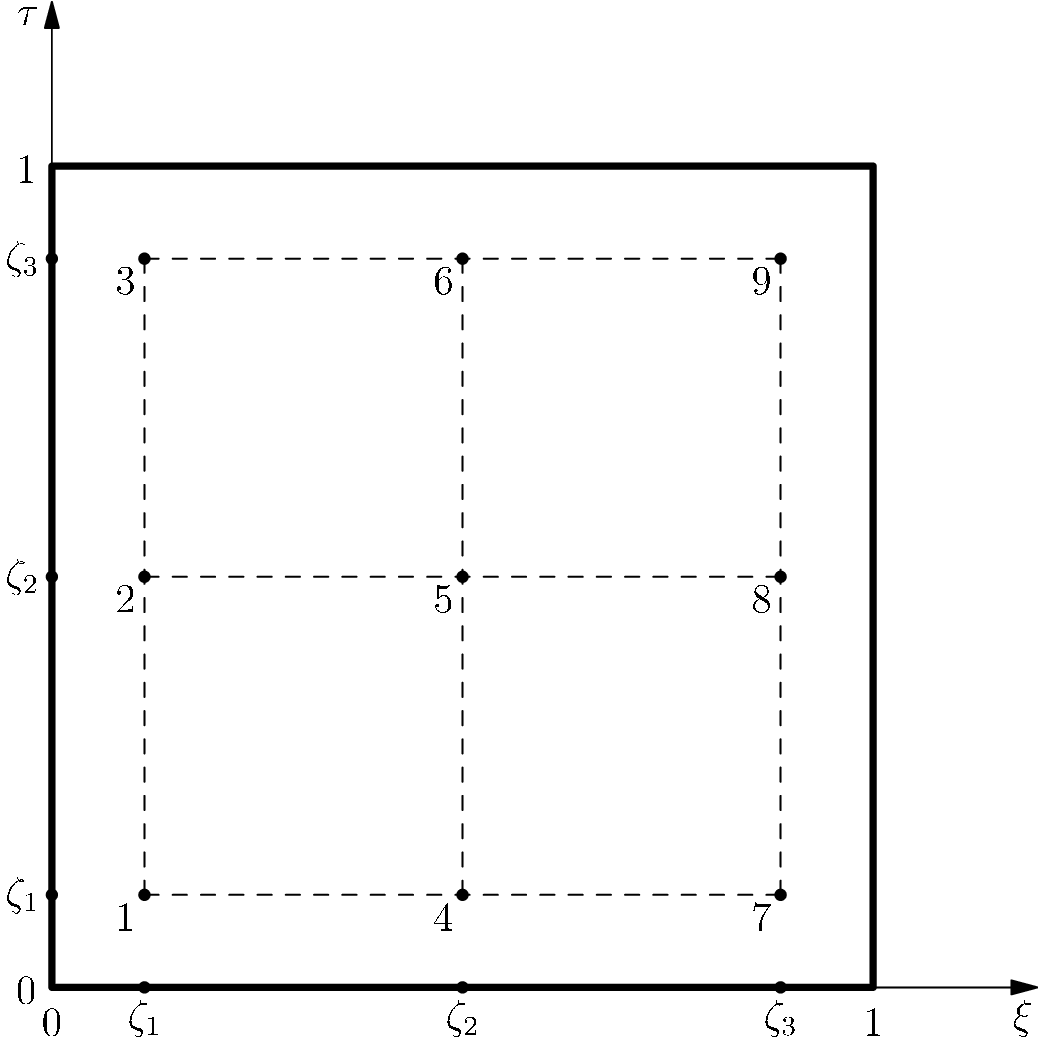}  
\end{tabular}
\caption{Sketch of a third order isoparametric space-time element. Left: physical space-time element. Right: reference space-time element. The interpolation nodes for 
the numerical solution and for the mapping, given by the tensor-product Gauss-Legendre quadrature points, are numbered from 1 to 9. The initial location for the spatial
Gauss-Legendre nodes $\hat x_{0,m}$ is also highlighted. }
\label{fig.isoparam}
\end{center}
\end{figure}

To get an element-local weak formulation of the PDE on a \textit{moving} space-time control volume $\mathcal{C}_i = [x_{i-\halb}(t);x_{i+\halb}(t)]\times[t^n;t^{n+1}]$, the 
governing PDE \eqref{eqn.pde} is mapped to the reference space-time element $\mathcal{C}_E=[0;1]^2$ using an \textit{isoparametric} mapping, hence the mapping of the geometry is 
approximated with the same space-time basis functions $\theta_m$ that are also used to approximate the discrete solution $\q_h(x,t)$. In this paper the $\theta_m$ are chosen to be 
the \textit{Lagrange} interpolation polynomials of degree $M$ that pass through the tensor-product \textit{Gauss-Legendre} quadrature points on the reference element $\mathcal{C}_E$, 
see \cite{stroud} for details on multidimensional quadrature. At the aid of the space-time basis functions $\theta_k$, the mapping of $x$ and $t$ onto $\xi$ and $\tau$ reads  
\begin{equation} 
\label{eqn.iso.map} 
 x_h = x(\xi,\tau) = \hat x_m \theta_m(\xi,\tau), \qquad  
 t_h = t(\xi,\tau) = \hat t_m \theta_m(\xi,\tau).  
\end{equation} 
Here, $\theta_m = \theta_m(\xi,\tau)$ and the coefficients $\hat x_m$ and $\hat t_m$ denote the \textit{nodal} coordinates in physical space and time and $\xi$ and $\tau$ 
are the reference coordinates. A sketch of this isoparametric mapping is depicted in Fig. \ref{fig.isoparam}. 
For the time coordinate one gets the following simple mapping, since the time coordinates are the same for each spatial node at each time level: 
\begin{equation}
\label{eqn.time.map} 
 t_h = t^n + \Delta t_i^n \tau.  
\end{equation} 
This reduces the Jacobian of the space-time mapping $(\xi,\tau) \to (x,t)$ and its inverse to  
\begin{equation}
 J = \left( \begin{array}{cc} x_\xi & x_\tau \\ t_\xi & t_\tau \end{array} \right) = \left( \begin{array}{cc} x_\xi & x_\tau \\ 0 & \Delta t_i^n \end{array} \right), \qquad    
J^{-1} = \left( \begin{array}{cc} \xi_x & \xi_t \\ \tau_x & \tau_t  \end{array} \right)  
       = \left( \begin{array}{cc} \frac{1}{x_\xi} & - \frac{1}{\Delta t_i^n} \frac{x_\tau}{x_\xi} \\ 0 & \frac{1}{\Delta t_i^n} \end{array} \right). 
\label{eqn.invjac} 
\end{equation}  
 
Using the chain rule PDE \eqref{eqn.pde} is rewritten on the reference element $\mathcal{C}_E$ as 
\begin{equation}
  \frac{\d \Q}{\d \tau} - \frac{x_\tau}{x_\xi} \frac{\d \Q}{\d \xi}  + \frac{\Delta t_i^n}{x_\xi} \frac{\d \F}{\d \xi} = 0,   
\label{eqn.pde.ref}  
\end{equation} 
where the second term corresponds to the Lagrangian part of the flux caused by the motion of the mesh. 
In the framework of isoparametric finite elements, the \textit{discrete} solution and the flux of PDE \eqref{eqn.pde.ref} are approximated with the same basis 
functions $\theta_m$ used for the mapping \eqref{eqn.iso.map}, i.e. 
\begin{equation}
\label{eqn.qh} 
\q_h = \q_h(\xi,\tau) = \theta_m(\xi,\tau) \hat \q_m, \qquad   
\f_h = \f_h(\xi,\tau) = \theta_m(\xi,\tau) \hat \f_m. 
\end{equation} 
For a nodal basis, the degrees of freedom of the interpolation of the nonlinear flux are simply computed pointwise as 
\begin{equation}
\label{eqn.nodalapprox} 
 \hat \F_m = \F(\hat \q_m). 
\end{equation} 
Throughout the paper, the Einstein summation convention is used. 
To ease notation, the following operators on $\mathcal{C}_E$ are introduced: 
\begin{equation}
\label{eqn.op.ref} 
  \left[f,g\right]^\tau  = \int \limits_0^1 f(\xi,\tau) g(\xi,\tau) \, d\xi, \qquad  \textnormal{ and } \qquad 
  \left<f,g\right>       = \int \limits_{0^+}^1 \int \limits_0^1 f(\xi,\tau) g(\xi,\tau) \, d\xi \, d\tau.  
\end{equation}  
Multiplication of Eqn. \eqref{eqn.pde.ref} with space-time test functions $\theta_k$, integration over $\mathcal{C}_E$ and integration of the first term
by parts in time yields 
\begin{equation}
\label{eqn.weak1} 
\left[ \theta_k, \q_h \right]^{1} - \left[ \theta_k, \w_h \right]^{0} - \left< \frac{\partial}{\partial \tau} \theta_k,  \q_h \right>  + 
\left< \theta_k, \frac{\Delta t_i^n}{x_\xi}  \frac{\partial}{\partial \xi} \F_h(\q_h) - \frac{x_\tau}{x_\xi} \frac{\d \q_h}{\d \xi} \right> = 0.
\end{equation}
Here, the initial condition given by the WENO reconstruction polynomial $\w_h(x,t^n) = \psi_m(\xi) \hat{\w}_m^n$ at time $t^n$ (see the next sub section) has been introduced in a \textit{weak form}. 

With the definitions for the WENO reconstruction polynomial and the discrete space-time solution \eqref{eqn.qh},  
one obtains the following element-local nonlinear algebraic equation system:  
\begin{equation}
\label{eqn.weak3} 
 K^1_{km} \hat \q_m + K^{\xi_x}_{km} \hat \F_m - K^{\xi_t}_{km} \hat \q_m = F^0_{km} \hat \w_m^n,  
\end{equation} 
with   
\begin{equation}
  K^1_{km} = \left( \left< \theta_k,  \frac{\partial}{\partial \tau} \theta_m \right> + \left[ \theta_k, \theta_m \right]^0 \right),  
\end{equation}  
\begin{equation}
 K^{\xi_x}_{km} = \left< \theta_k,  \frac{\Delta t_i^n}{x_\xi} \frac{\partial \theta_m}{\partial \xi} \right>, \qquad 
 K^{\xi_t}_{km} = \left< \theta_k,  \frac{x_\tau}{x_\xi}   \frac{\partial \theta_m}{\partial \xi} \right>, 
\end{equation}  
and 
\begin{equation}
F^0_{km} = \left[\theta_k, \psi_m \right]^0, \qquad 
M_{km} = \left< \theta_k, \theta_m  \right>. 
\end{equation} 
The element-local nonlinear algebraic systems \eqref{eqn.weak3} can be easily solved using the following iterative method, see \cite{DumbserZanotti,HidalgoDumbser}: 
\begin{equation}
\label{eqn.weak.iter} 
 K^1_{km} \hat \q_m^{l+1} + K^{\xi_x}_{km} \hat \F_m^{l} - K^{\xi_t}_{km} \hat \q_m^{l} = F^0_{km} \hat \w_m^n.  
\end{equation} 
For an efficient strategy to get the initial guess $\hat \q^0_m$ see \cite{HidalgoDumbser}. 
The equation that determines the location of the spatial coordinates $\hat x_m$ of the space-time element is the ODE 
\begin{equation}
\label{eqn.mesh.motion} 
  \frac{dx}{dt} = V(\Q(x,t)),  
\end{equation} 
where $V(\Q(x,t))$ is the local mesh velocity. For the local mesh velocity one can use again the nodal ansatz 
\begin{equation}
  V_h = V(\q_h(x,t)) = \theta_m(x,t) \hat v_m, \qquad \textnormal{ with } \qquad \hat v_m = V(\hat{\q}_m). 
\end{equation}
The initial distribution of the spatial Gauss-Legendre quadrature points at time $t_i^n$ is given by
\begin{equation}
  \hat x_{0,m} = x_{i-\halb}^n + \Delta x_i^n \zeta_m, 
\end{equation} 
where the $\zeta_m$ are the quadrature points on the unit interval $[0;1]$ and the spatial Lagrange interpolation polynomials 
passing through these points are denoted by $\phi_m$. 
A discrete version of the ODE \eqref{eqn.mesh.motion} can then be obtained using again the local space-time DG method, 
see \cite{ADERNSE}: 
\begin{equation}
\label{eqn.weak.x} 
  \left(  [\theta_k, \theta_m]^1 - \left< \frac{\partial}{\partial \tau}  \theta_k, \theta_m \right>  \right) \hat x^{l+1}_m = [\theta_k, \phi_m]^0 \hat x_{0,m} + 
   \Delta t_i^n \left< \theta_k, \theta_m \right> \hat v^l_m.  
\end{equation} 
The weak formulation for the spatial coordinates \eqref{eqn.weak.x} is iterated \textit{together} with the weak formulation for the solution \eqref{eqn.weak.iter} until 
convergence is reached. The temporal coordinates $\hat t_m$ are \textit{fixed} and are given by the Gauss-Legendre points $\zeta_m$ in time and relation \eqref{eqn.time.map}.  
The space-time polynomials $\q_h(x,t)$ are computed for each element in the computational domain and are then used as arguments for the numerical flux in Eqn. \eqref{eqn.numflux}.

\subsection{High Order WENO Reconstruction for LTS} 
\label{sec.weno} 

In this paper the \textit{polynomial} WENO reconstruction algorithm proposed in \cite{DumbserKaeser06b,DumbserKaeser07,DumbserEnauxToro} is used. 
Its output are reconstruction polynomials and not point values, as in the original optimal WENO scheme of Jiang and Shu \cite{shu_efficient_weno}. 
The details can be found in the above-mentioned references, hence only a brief summary of the algorithm is given here together with the necessary 
modifications in order to support LTS. For the sake of simplicity, componentwise reconstruction in conservative variables is assumed. 
For reconstruction in characteristic variables see \cite{eno,shu_efficient_weno}.  
A reconstruction polynomial $\w_h(x,t_i^n)$ of degree $M$ is obtained componentwise by requiring \textit{integral conservation} on a stencil 
\begin{equation}
\label{eqn.stencildef}  
\mathcal{S}_i^{s} = \bigcup \limits_{j=i-l}^{i+r} T_j^n   
\end{equation}
with spatial extension $l$ and $r$ to the left and right, respectively. For odd order schemes there is one central stencil ($s=1$), with $l=r=M/2$, 
while for even order schemes there are two central stencils with $l=$floor$(M/2)+1$ and $r=$floor$(M/2)$ for the first central stencil ($s=0$) and 
$l=$floor$(M/2)$, $r=$floor$(M/2)+1$ for the second one ($s=1$). All schemes have one fully left-sided stencil ($s=2$) with $l=M$ and $r=0$ and 
one fully right-sided stencil ($s=3$) with $l=0$ and $r=M$. The reconstruction polynomial for each stencil is written in terms of some spatial basis functions 
$\psi_m(\xi)$ as  
\begin{equation}
\label{eqn.recpolydef} 
 \w^s_h(x,t^n) = \sum \limits_{m=0}^M \psi_m(\xi) \hat \w^s_m := \psi_m(\xi) \hat \w^{n,s}_m,   
\end{equation}
with the mapping  
\begin{equation}
 x = x_{i-\halb}^n + \Delta x_i^n \xi. 
\end{equation} 
For the reconstruction basis functions $\psi_m(\xi)$ one can either use Legendre polynomials rescaled to the unit interval $[0;1]$, or, a nodal basis based
on the Lagrangian interpolation polynomials passing through the Gauss-Legendre quadrature points on the unit interval. To obtain the high order 
reconstruction polynomial one usually requires integral conservation on all elements of the stencil as follows: 
\begin{equation}
 \frac{1}{\Delta x^n_j} \int \limits_{x^n_{j-\halb}}^{x^n_{j+\halb}} \w^s_h(x,t_i^n) dx = \Q^n_j, \qquad \forall T_j^n \in \mathcal{S}_i^{s}.     
 \label{eqn.receqn1} 
\end{equation}
However, in the context of finite volume schemes with time-accurate local time stepping, each cell average $\Q^n_j$ is usually defined at a \textit{different} 
local time $t_j^n \neq t_i^n$, i.e. Eqn. \eqref{eqn.receqn1} can in general \textit{not} be used for reconstruction, since the cell averages needed 
for reconstruction are in general not available at the time $t_i^n$ when they are needed. Also the geometry of the stencil is in general not available at
the time $t_i^n$ required for reconstruction. Instead, one has to use the following set of reconstruction equations, which uses the local predictor 
solution $\q_h$ in all elements $T_j^n$ with $j \neq i$ in order to \textit{predict} the value of the cell averages and the geometry at the required
reconstruction time $t_i^n$. These estimated values are distinguished from the real cell averages and interface positions by the use of the tilde symbol:     
\begin{equation}
 \frac{1}{\Delta \tilde x^n_j} \int \limits_{\tilde x_{j-\halb}(t_i^n)}^{\tilde x_{j+\halb}(t_i^n)} \w^s_h(x,t_i^n) dx = \tilde{\Q}_j(t_i^n), \qquad \forall T_j^n \in \mathcal{S}_i^{s},     
 \label{eqn.receqn2} 
\end{equation}
with $\Delta \tilde x^n_j = \Delta \tilde x_j(t_i^n) = \tilde x_{j+\halb}(t_i^n) - \tilde x_{j-\halb}(t_i^n)$ and  
\begin{equation}
  \tilde{\Q}_j(t_i^n) = 
	\left\{ \begin{array}{lll} 
	\Q^n_i, &  \textnormal{ if }  & i=j, \\ 
	\frac{1}{\Delta \tilde x^n_j} \displaystyle \int \limits_{\tilde x_{j-\halb}(t_i^n)}^{\tilde x_{j+\halb}(t_i^n)} \hspace{-3mm} \q_h(x,t_i^n) dx, &  \textnormal{ if } & i \neq j. 
	\end{array} \right. 
  \label{eqn.virtual.qh} 
\end{equation} 
The geometry at the reconstruction time $t_i^n$ is estimated by moving the interfaces \textit{virtually}, just for the purpose of reconstruction, as follows: 
\begin{equation} 
\tilde x_{j+\halb}(t_i^n) = x_{j+\halb}(t_{j+\halb}^n) + \int \limits_{t_{j+\halb}^{n}}^{t_i^{n}} V_{j+\halb}(t) dt, 
  \label{eqn.virtual.x} 
\end{equation} 
where $V_{j+\halb}(t) = \frac{1}{2}\left( V(\q_h(x_{j+\halb}^-(t),t))+V(\q_h(x_{j+\halb}^+(t),t))\right)$ 
is again easily computed using the local predictor $\q_h$. Now, the reconstruction equations \eqref{eqn.receqn2} are defined and the unknown coefficients 
$\w^{n,s}_m$ can be obtained, using the \textit{virtual} interface positions $\tilde x_{j+\halb}(t_i^n)$ and the virtual cell averages $\tilde{\Q}^n_j$ 
at the reconstruction time $t_i^n$, computed from the element-local space-time predictor solution $\q_h(x,t)$. As the reader can easily appreciate, it is indeed the 
use of the predictor-corrector philosophy which allows a simple reconstruction procedure with time-accurate local time-stepping, since the local space-time predictor 
$\q_h$ provides an accurate interpolation of the data and the geometry at any desired position in space and time. From \eqref{eqn.virtual.x} and \eqref{eqn.virtual.qh} 
it becomes clear that an element can only apply the reconstruction operator if the following condition is satisfied: 
\begin{equation}
\label{eqn.stencil.condition}
 \max( t_j^n ) \leq t_i^n \leq \min( t_j^{n+1} ), \qquad \forall T_j^n \in \mathcal{S}_i^s, 
\end{equation} 
i.e. the reconstruction time $t_i^n$ must be contained within all local time intervals $[t_j^n;t_j^{n+1}]$ of all stencil elements so that the local predictor solution can be evaluated at a valid relative time $0 \leq \tau \leq 1$. The predictor is not valid for $\tau<0$ (backward extrapolation in time) or $\tau>1$ 
(forward extrapolation in time). Note that for better than second order schemes condition \eqref{eqn.stencil.condition} is more restrictive than the simple 
update criterion \eqref{eqn.update.criterion}, hence from third order methods onward one has to obey condition \eqref{eqn.stencil.condition}.   

With the oscillation indicators $\sigma_s$ 
\begin{equation}
\sigma_s = \Sigma_{lm} \hat \w^{n,s}_l \hat \w^{n,s}_m, \qquad 
\Sigma_{lm} = \sum \limits_{\alpha=1}^M \int \limits_0^1 \frac{\partial^\alpha \psi_l(\xi)}{\partial \xi^\alpha} \cdot \frac{\partial^\alpha \psi_m(\xi)}{\partial \xi^\alpha} d\xi,   
\end{equation} 
the nonlinear weights $\omega_s$ are defined by
\begin{equation}
\tilde{\omega}_s = \frac{\lambda_s}{\left(\sigma_s + \epsilon \right)^r}, \qquad 
\omega_s = \frac{\tilde{\omega}_s}{\sum_q \tilde{\omega}_q},  
\end{equation} 
where we use $\epsilon=10^{-14}$, $r=8$, $\lambda_s=1$ for the one-sided stencils and $\lambda=10^5$ for the central stencils, according to \cite{DumbserEnauxToro,DumbserKaeser06b}. 
The final nonlinear WENO reconstruction polynomial and its coefficients are then given by 
\begin{equation}
\label{eqn.weno} 
 \w_h(x,t^n) = \psi_m(\xi) \hat \w^{n}_m, \qquad \textnormal{ with } \qquad  
 \hat \w^{n}_m = \sum_s \omega_s \hat \w^{n,s}_m.   
\end{equation}   

For the computation of the time step $\Delta t_i$, each element obeys the following 
\textit{local CFL condition}
\begin{equation}
\label{eqn.localcfl} 
  \Delta t_i^n = \textnormal{CFL} \cdot \min_{j \in \mathcal{N}_i} \left(  \frac{\Delta \tilde{x}_j(t_i^n)}{ \left| \lambda_{\max}( \tilde{\Q}_j(t_i^n) )\right| } \right), 
	\qquad \textnormal{with} \qquad \textnormal{CFL} \leq 1.  
\end{equation} 

This closes the description of the high order reconstruction operator in the context of LTS and thus completes the general presentation of the algorithm, which consists in the three steps
i) reconstruction, ii) data evolution and iii) cell update. 

\subsection{Description of the Local Time-Stepping Algorithm} 
\label{sec.lts} 

\begin{figure}[!t] 
\begin{center}
\begin{tabular}{lr}
\includegraphics[width=0.45\textwidth]{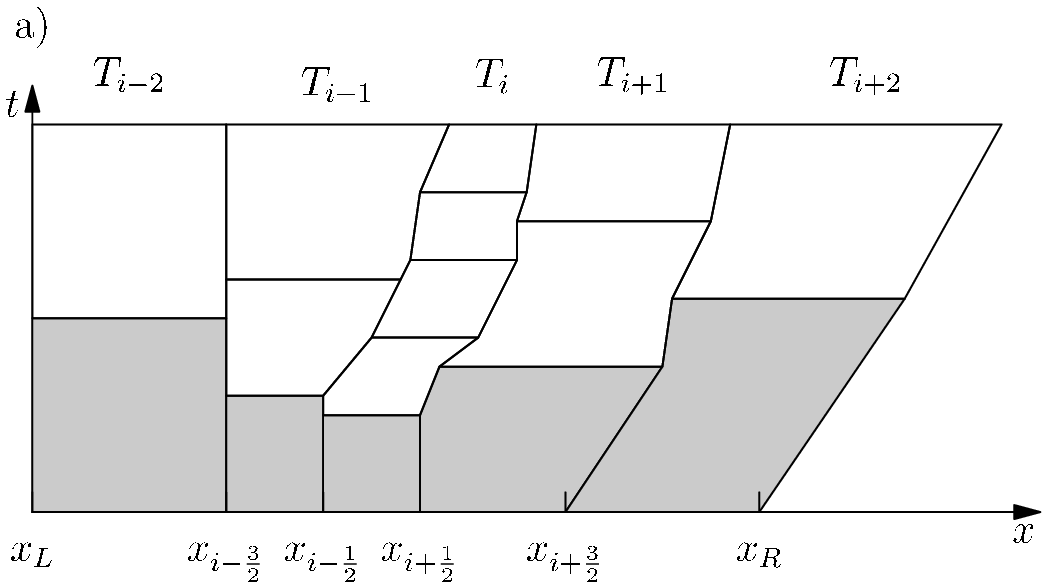}  & 
\includegraphics[width=0.45\textwidth]{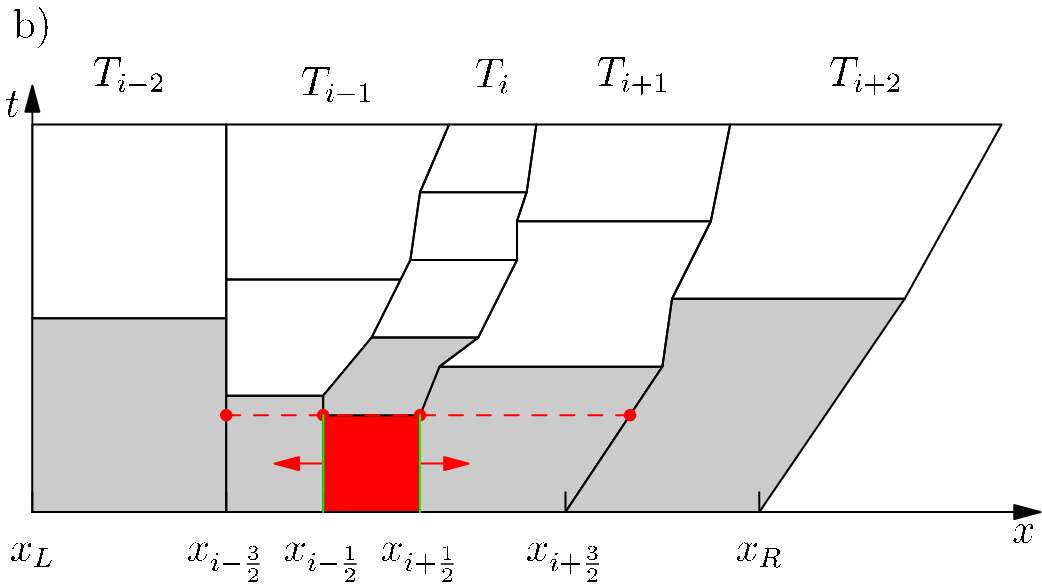}  \\
\includegraphics[width=0.45\textwidth]{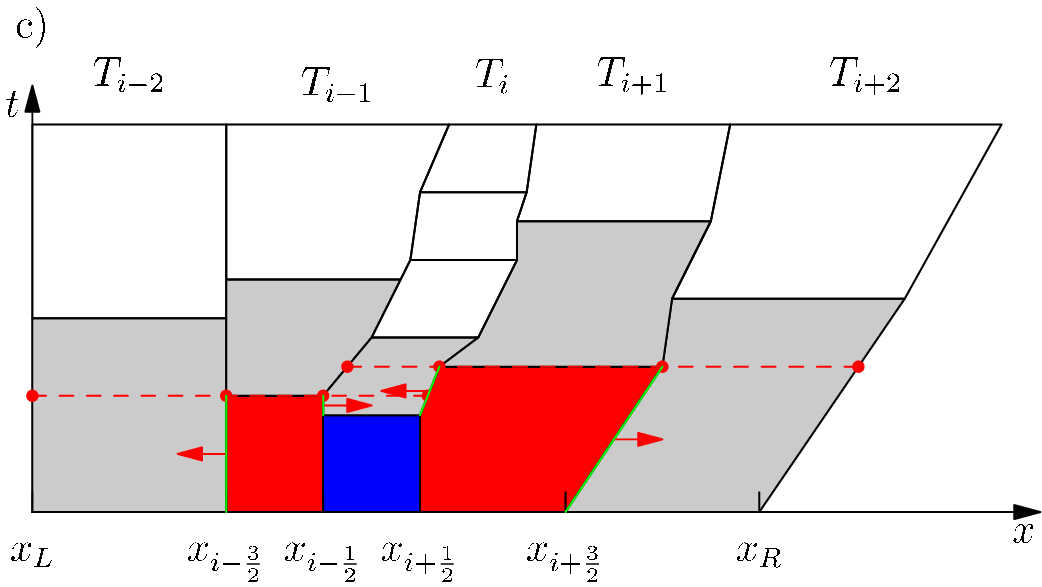}  & 
\includegraphics[width=0.45\textwidth]{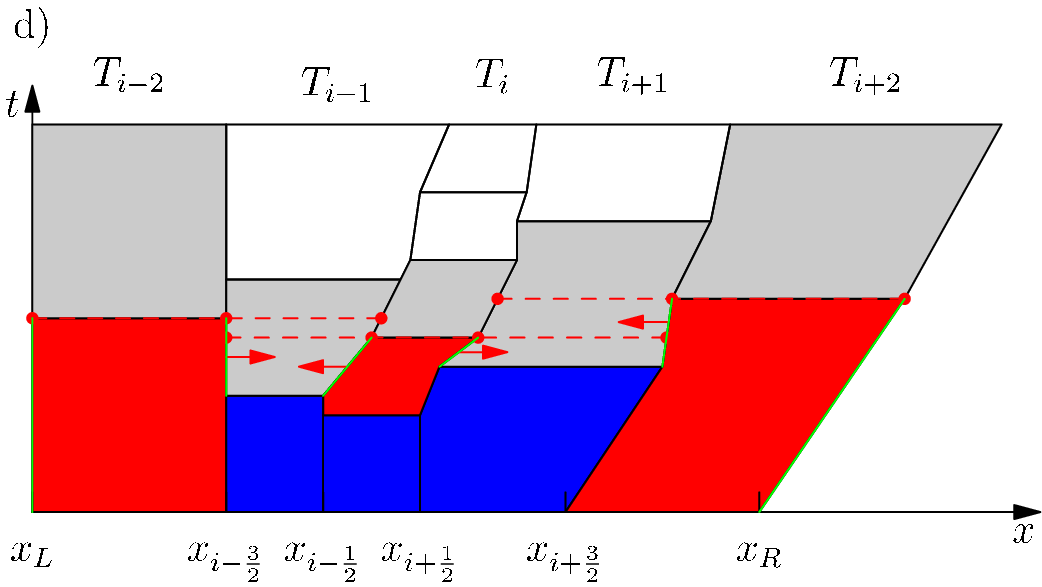}   
\end{tabular}
\caption{Sketch of the ALE-LTS procedure and the generated space-time mesh. Red space-time cells are currently being updated, blue space-time cells have already been updated at least once and space-time cells with an available predictor solution $\q_h$ are highlighted 
in grey. The flux time intervals are depicted in green, the red arrows denote saving of 
flux contributions in memory variables of neighbor cells and the horizontal dashed 
red lines indicate reconstruction using virtual data obtained by interpolation from 
the space-time predictor $\q_h$. Note the hanging nodes in time 
produced by the ALE-LTS algorithm on the element edges.} 
\label{fig.lts}
\end{center}
\end{figure}

The workflow of an LTS algorithm is best illustrated on a simple example, as the one sketched in Figure \ref{fig.lts}. For other examples, see 
\cite{dumbserkaeser06d,TaubeMaxwell,stedg1}. 
Note that in an LTS algorithm the order in which the elements are updated depends on the update criterion \eqref{eqn.stencil.condition}. For that reason, one cannot speak of 
\textit{timesteps} any more, but the time marching of the LTS code is organized by \textit{cycles}. In each cycle, the scheme runs over \textit{all} elements of the domain and skips 
those which do not satisfy the update criterion \eqref{eqn.stencil.condition}. This verification is very fast and does not add much computational overhead to the scheme. For the 
example used here, we suppose $M=1$, hence a second order scheme in space and time is used, for which the criteria \eqref{eqn.update.criterion} and \eqref{eqn.stencil.condition} 
are the same. In this example, suppose that elements $T_{i-2}$ and $T_{i+2}$ are on the
boundary of the domain, hence they have to fulfill the update criterion only with respect
to neighbors within the computational domain.  

At the beginning of the time marching, all elements are defined at the common initial 
time $t=0$. The cell averages are defined by the initial condition and the \textit{first} 
WENO reconstruction can be done for all elements directly based on the cell averages 
$\Q_j$, without the use of the virtual geometry or the virtual cell averages 
$\tilde{\Q}_j$. For each element an element-local time step is defined according to the 
local CFL condition \eqref{eqn.localcfl}. Note that for the computation of the first time 
step, the real geometry and the real cell averages can be used in \eqref{eqn.localcfl} 
instead of the virtual ones. After reconstruction and local time step computation, each element can 
compute its element-local predictor solution $\q_h$. We are now in the situation depicted in 
Figure \ref{fig.lts} a). All cells are at the same initial time $t=0$, have different 
time steps and have computed their predictor $\q_h(x,t)$, which is valid through one local time 
step. The intervals with a valid space-time predictor solution $\q_h(x,t)$ are highlighted 
in grey in Figure \ref{fig.lts}. 
In the first cycle, only element $T_i$ satisfies the update criterion \eqref{eqn.stencil.condition}, hence it is updated according to the finite volume scheme \eqref{eqn.fv}, since on both element interfaces $x_{i \pm \halb}$ the necessary data for 
flux calculation is available under the form of the predictor $\q_h$. The updated cell is
highlighted in red in Figure \ref{fig.lts} b) and the corresponding flux intervals in 
time according to Eqn. \eqref{eqn.flux.timeint} are highlighted in green. Note that
the geometry is updated only for the nodes attached to the element which is updated,
hence $x_{i-\halb}$ and $x_{i+\halb}$ advance in time, while the other nodes are still
at their old time level. 
Once the fluxes for cell $T_i$ have been computed, they are added to the memory 
variables $\Q_{i\pm1}^M$ of the neighbor cells (red arrows in Figure \ref{fig.lts}) 
and the memory variable $\Q_i^M$ is reset to zero, according to \eqref{eqn.memvar.def}. 
Now, the local time of cell $T_i$ is \textit{advanced} by a local time step and  reconstruction is performed for cell $T_i$ at the \textit{new} local time $t_i^n$. 
Since reconstruction requires information from the neighbor cells $T_{i\pm1}$, which 
are still at the initial time $t=0$, the reconstruction operator uses the 
\textit{virtual data} defined in \eqref{eqn.virtual.qh} and a \textit{virtual geometry} 
according to \eqref{eqn.virtual.x}, which are both computed from the local space-time  predictor $\q_h(x,t)$. This reconstruction based on virtual data is 
highlighted by the dashed red horizontal lines in Figure \ref{fig.lts}. 
After reconstruction, cell $T_i$ can 
compute its predictor, highlighted again in grey in Figure \ref{fig.lts} b), and
the first cycle is complete, since no other cells satisfy the update criterion. 
The second cycle is depicted in Figure \ref{fig.lts} c). Now, cells $T_{i-1}$ and 
$T_{i+1}$ satisfy the update criterion and thus they are both updated in the 
\textit{same} cycle. The resulting flux time intervals according to 
Eqn. \eqref{eqn.flux.timeint} are highlighted again in green. One can note 
in Figure \ref{fig.lts} c) that the part of the fluxes which has already 
been computed previously by element $T_i$ and has been saved into the memory 
variables of $T_{i+1}$ and $T_{i-1}$ is not computed again, but only the missing
parts of the time intervals that are necessary to reach the future times $t_{i+1}^{n+1}$
and $t_{i-1}^{n+1}$, respectively. Again, the computed fluxes are added to the 
memory variables of the neighbor elements, hence $T_{i-1}$ adds the corresponding 
flux contributions to element $T_{i-2}$ on the left and to $T_i$ on the right, 
see the red arrows in Figure \ref{fig.lts} c). 
Likewise, element $T_{i+1}$ contributes to the memory variables of $T_i$ and
$T_{i+2}$. After resetting their memory variables, elements $T_{i-1}$ and $T_{i+1}$ 
perform the reconstruction using the virtual geometry and data provided by the 
space-time predictor $\q_h$ in the neighbor elements, see the dashed lines in 
Figure \ref{fig.lts}. After reconstruction, the predictor is computed again, 
see the grey shaded areas of elements $T_{i-1}$ and $T_{i+1}$. Since 
in the same cycle no other elements fulfill the update criterion, the cycle ends
and we have the situation depicted in Figure \ref{fig.lts} d). The three elements 
that have been updated so far are highlighted in blue, and the elements that satisfy
the update criterion within this cycle are $T_{i-2}$, $T_i$ and $T_{i+2}$. As
in the previous cycles, the flux time intervals are highlighted in green, the 
communication with the memory variables of the neighbor elements is indicated 
by the red arrows and the reconstruction based on the virtual data is 
highlighted by the red dashed horizontal lines. Recall also, that the space-time
mesh depicted in Figure \ref{fig.lts} is constructed \textit{dynamically} 
in each cycle according to the local CFL condition \eqref{eqn.localcfl} and 
the Lagrangian mesh motion. It is further important to note that only 
the spatial nodes attached to the elements that are updated (highlighted in red) 
are advanced in time. The other nodes remain at their old position at 
their old time. The resulting space-time mesh is non-conforming in time. 
As a direct consequence, hanging nodes in time appear, which are conveniently 
treated by the algorithm using memory variables and the fact that the flux 
time integrals are additive. 
The concept of memory variables has first been introduced for LTS methods 
in \cite{dumbserkaeser06d}. 

\section{Numerical Test problems} 
\label{sec.tests}

In this section, the numerical approach presented above is validated on a large set of one-dimensional Riemann problems for the Euler equations 
of compressible gas dynamics and the equations of ideal magnetohydrodynamics (MHD). Furthermore, a detailed numerical convergence study is presented 
for the MHD system for a smooth, non-trivial time-dependent test case. 

\subsection{Euler equations of compressible gas dynamics} 
\label{sec.euler}

For the first set of shock tube problems the classical Euler equations of compressible gasdynamics are considered. In one space
dimension they read
\begin{equation}
\label{eqn.Euler} 
  \frac{\partial}{\partial t} \left( \begin{array}{c} 
  \rho \\   \rho u \\   \rho E  \end{array} \right)   
  + \frac{\partial}{\partial x} \left( \begin{array}{c} 
  \rho u  \\ 
  \rho u^2 + p \\ 
  u (\rho E + p ) \\ 
  \end{array} \right) = 0, 
\end{equation} 
with 
\begin{equation}
  p = (\gamma - 1) (\rho E - \halb \rho u^2 ),
\end{equation} 
where $\rho$ is the mass density, $u$ is the velocity, $p$ is the gas pressure, $\rho E$ is the total energy density and $\gamma$ is the 
ratio of specific heats. 

\begin{figure}
\begin{center}
\begin{tabular}{lr}
\includegraphics[width=0.45\textwidth]{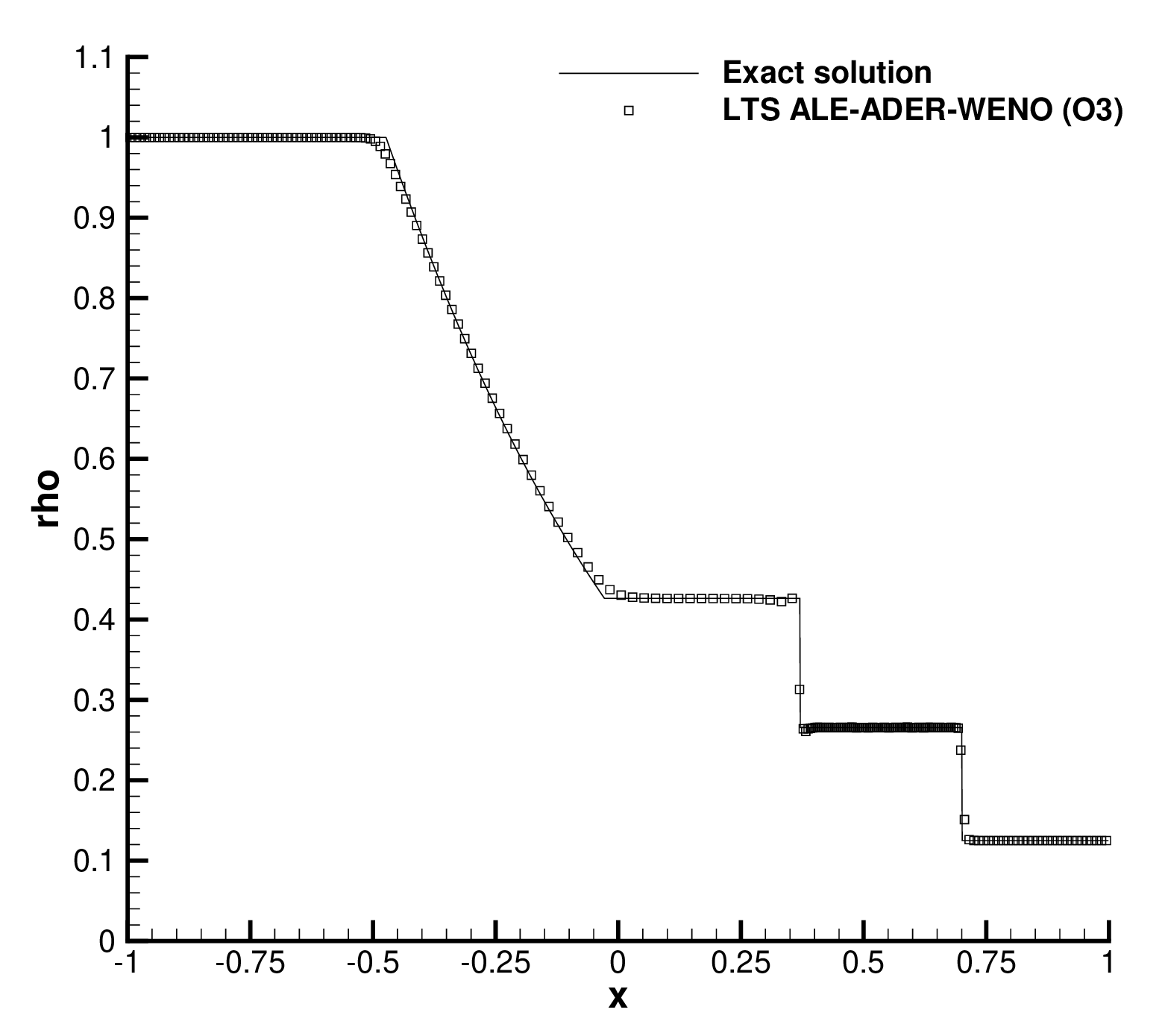} & 
\includegraphics[width=0.45\textwidth]{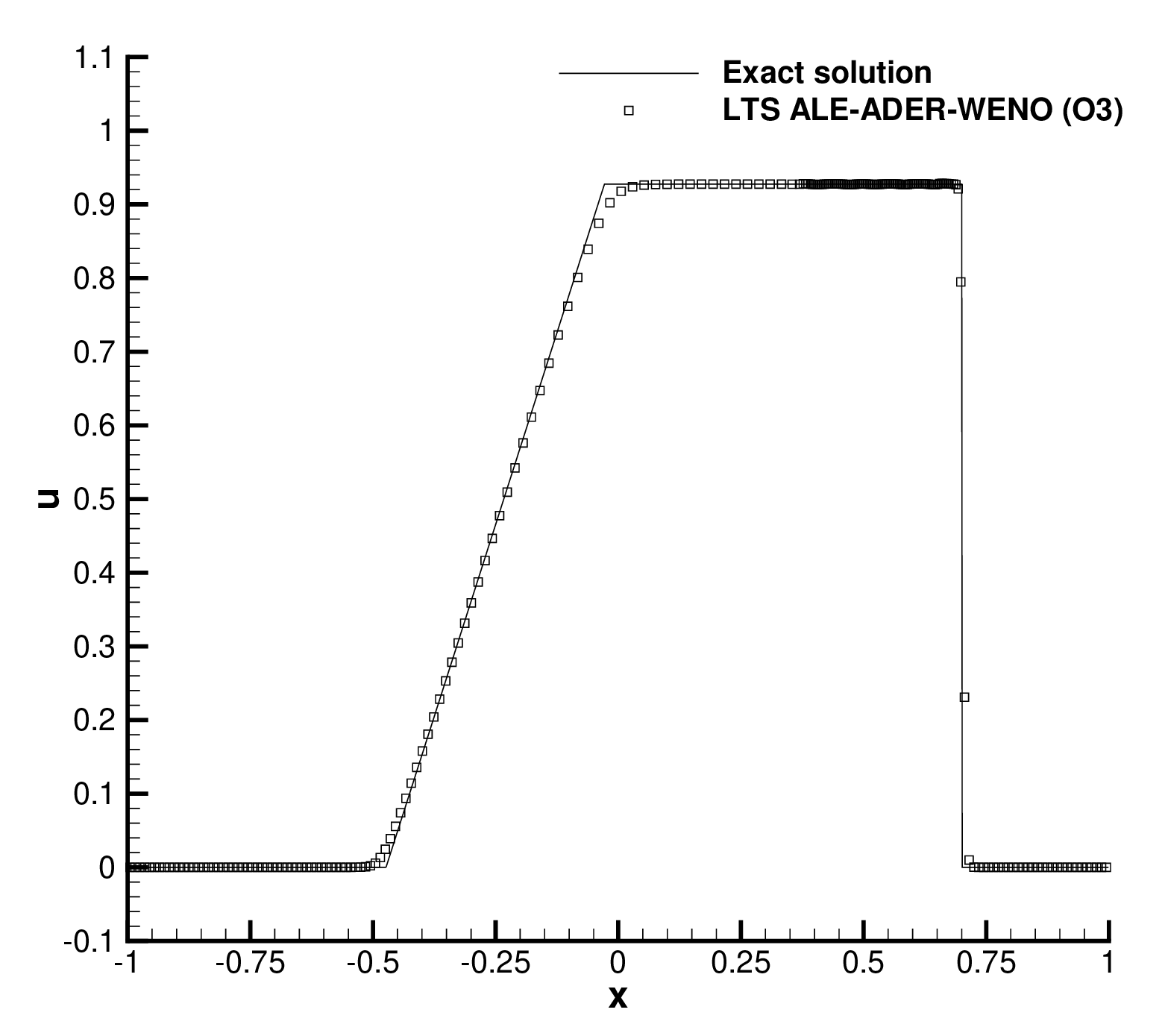}  \\
\multicolumn{2}{c}{\includegraphics[width=0.75\textwidth]{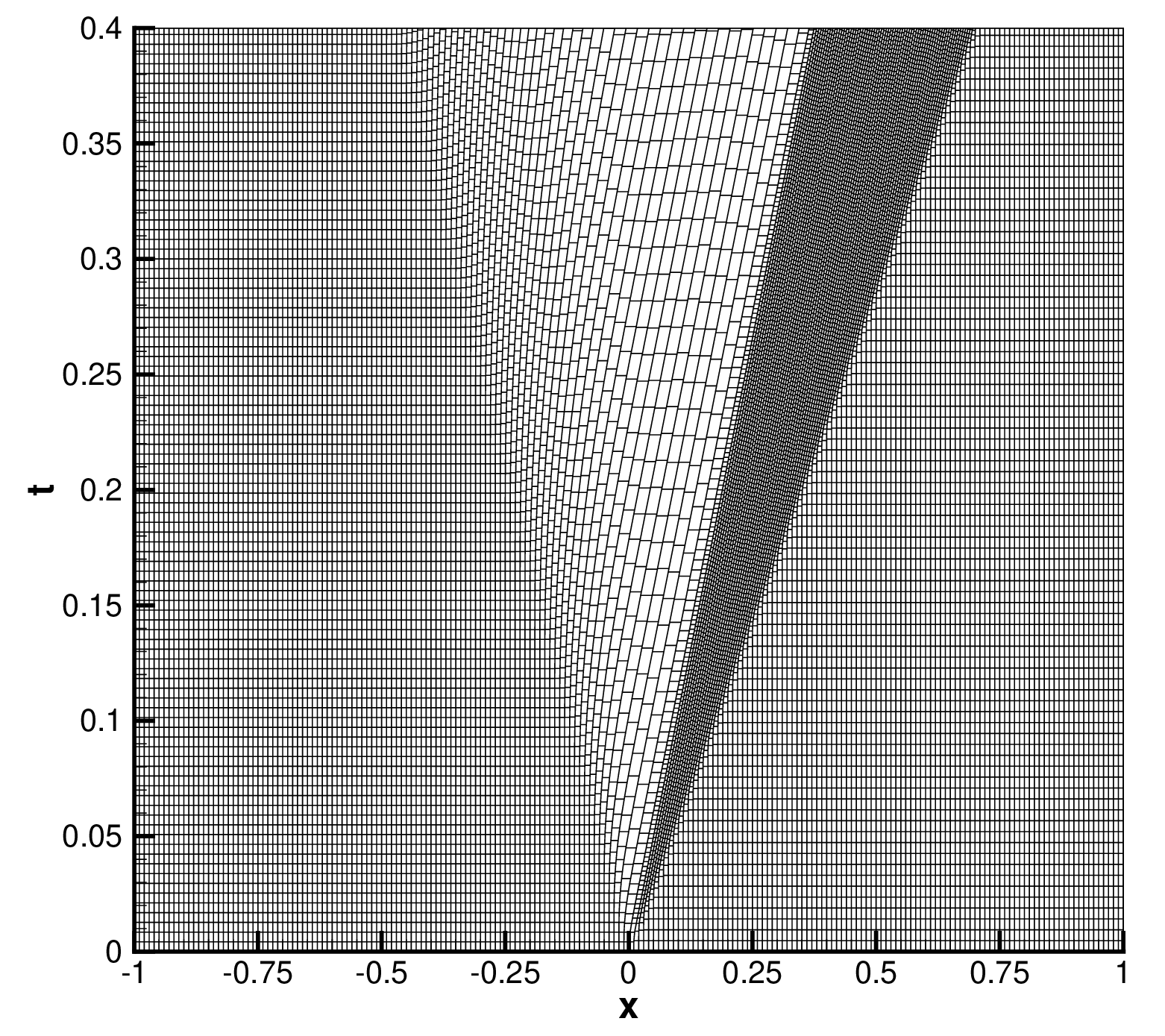}} 
\end{tabular}
\caption{Exact and numerical solution for the Sod shock tube problem. Density (top left), velocity (top right) and resulting space-time mesh of the third order Lagrangian ADER-WENO finite volume scheme with conservative and time-accurate local time stepping (bottom). }
\label{fig.sod}
\end{center}
\end{figure}

\begin{figure}
\begin{center}
\begin{tabular}{lr}
\includegraphics[width=0.45\textwidth]{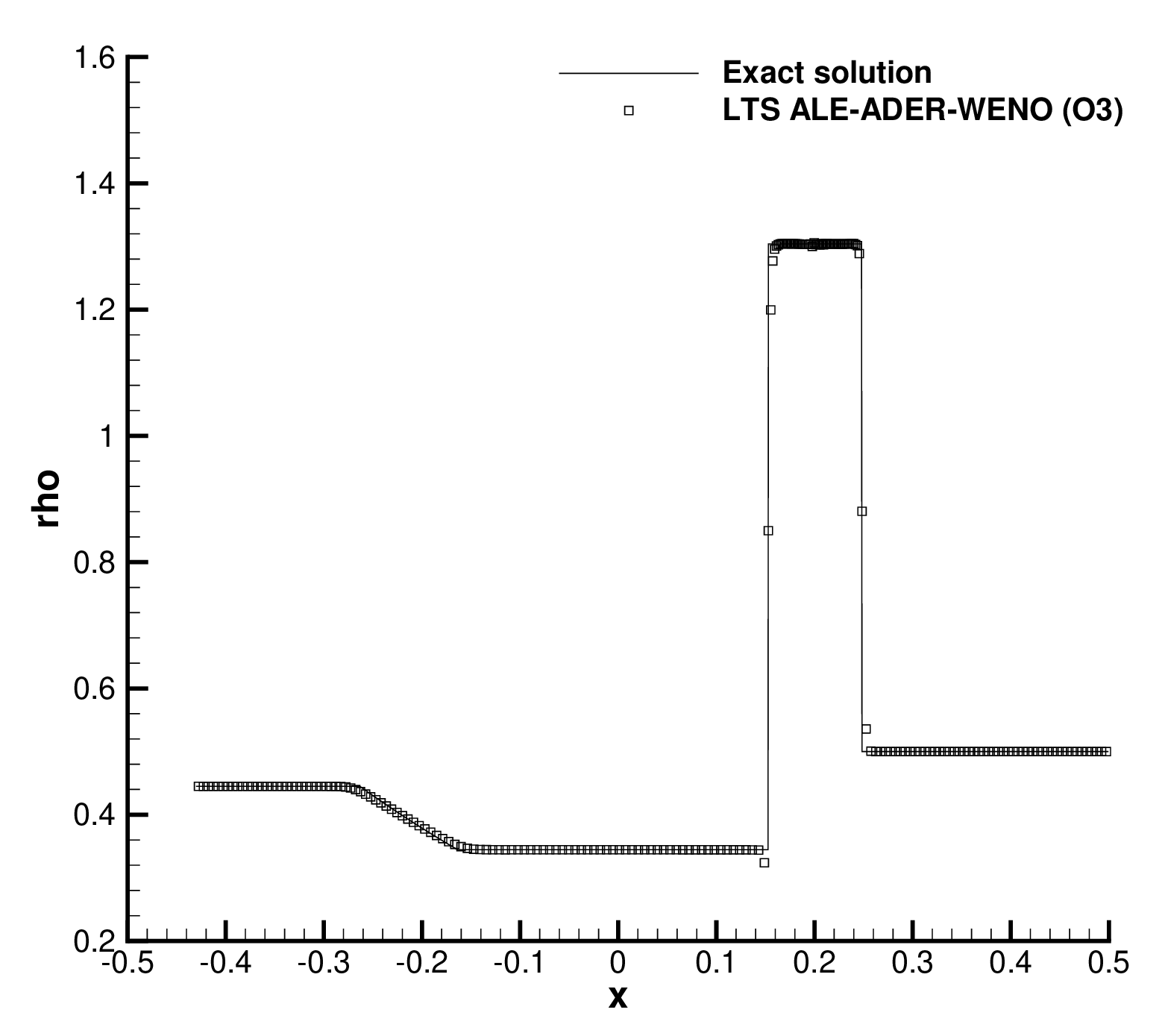} & 
\includegraphics[width=0.45\textwidth]{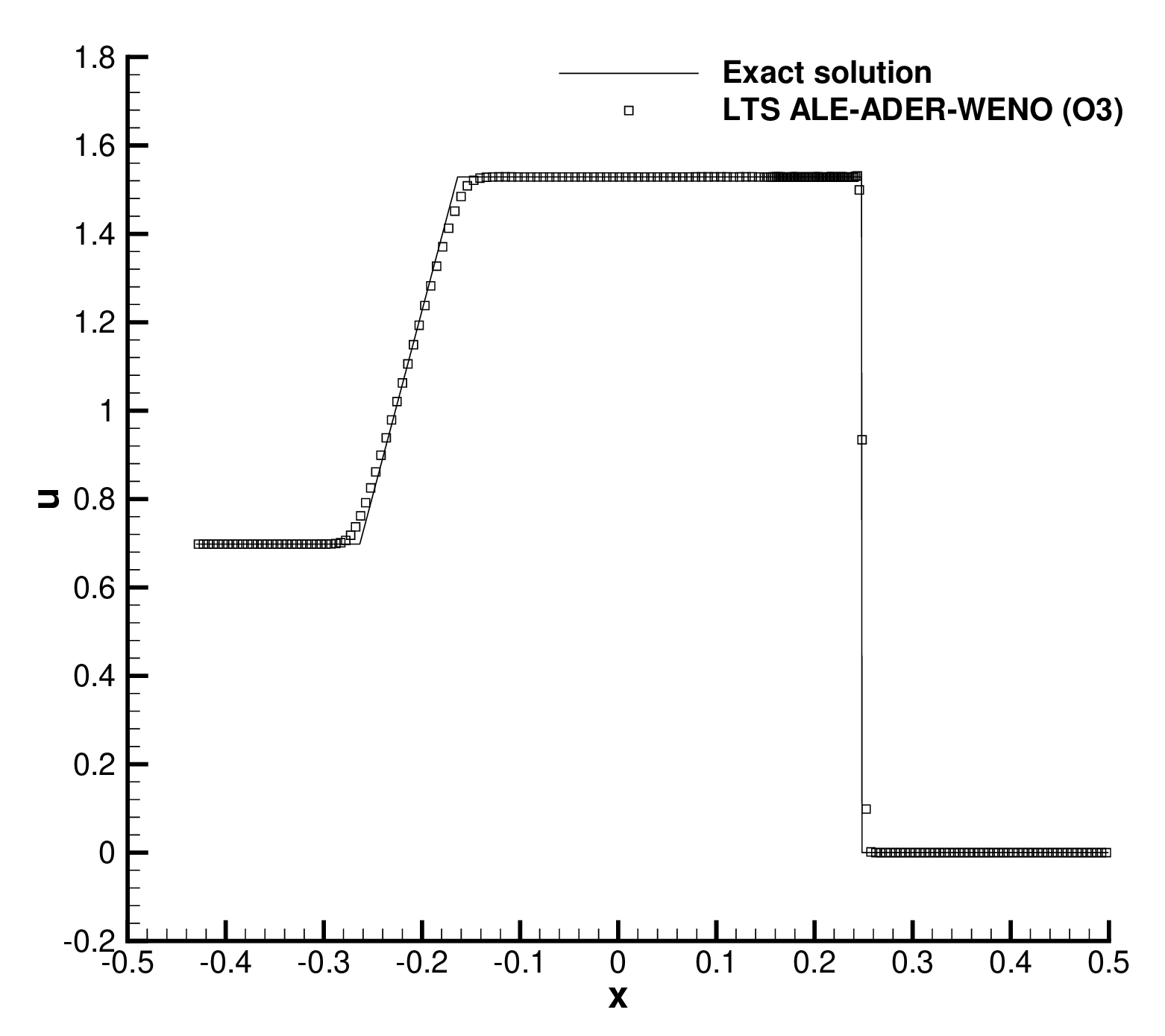}  \\
\multicolumn{2}{c}{\includegraphics[width=0.75\textwidth]{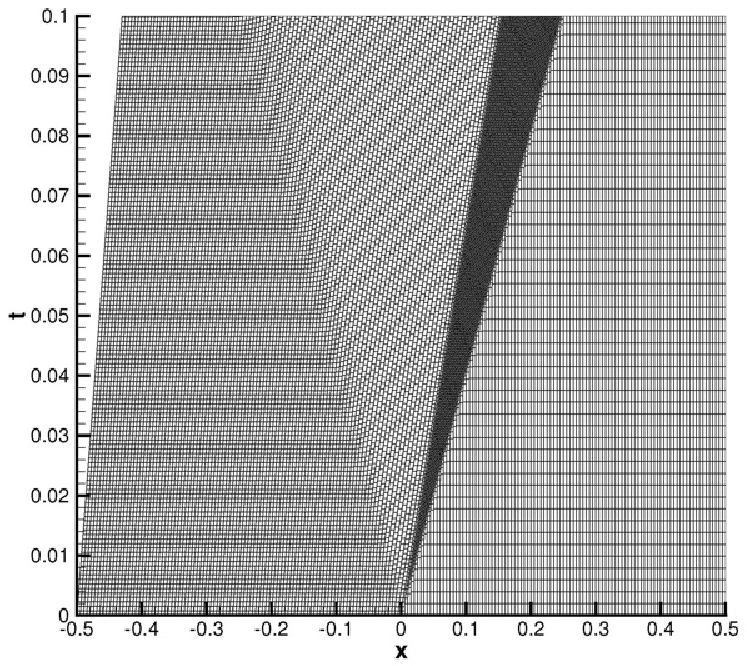}} 
\end{tabular}
\caption{Exact and numerical solution for the Lax shock tube problem. Density (top left), velocity (top right) and resulting space-time mesh of the third order Lagrangian 
ADER-WENO finite volume scheme with conservative and time-accurate local time stepping (bottom).}
\label{fig.lax}
\end{center}
\end{figure}

\begin{figure}
\begin{center}
\begin{tabular}{lr}
\includegraphics[width=0.45\textwidth]{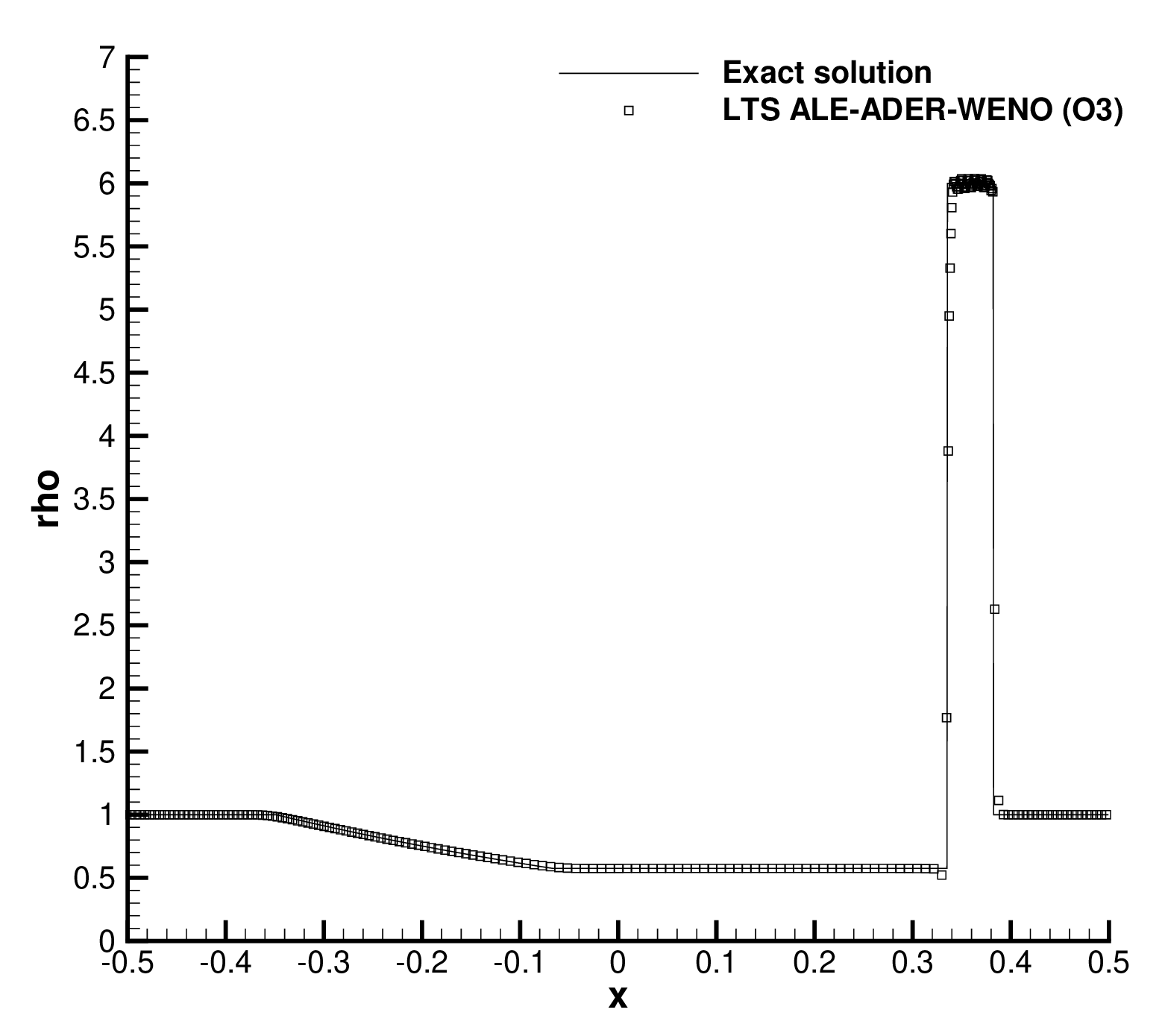} & 
\includegraphics[width=0.45\textwidth]{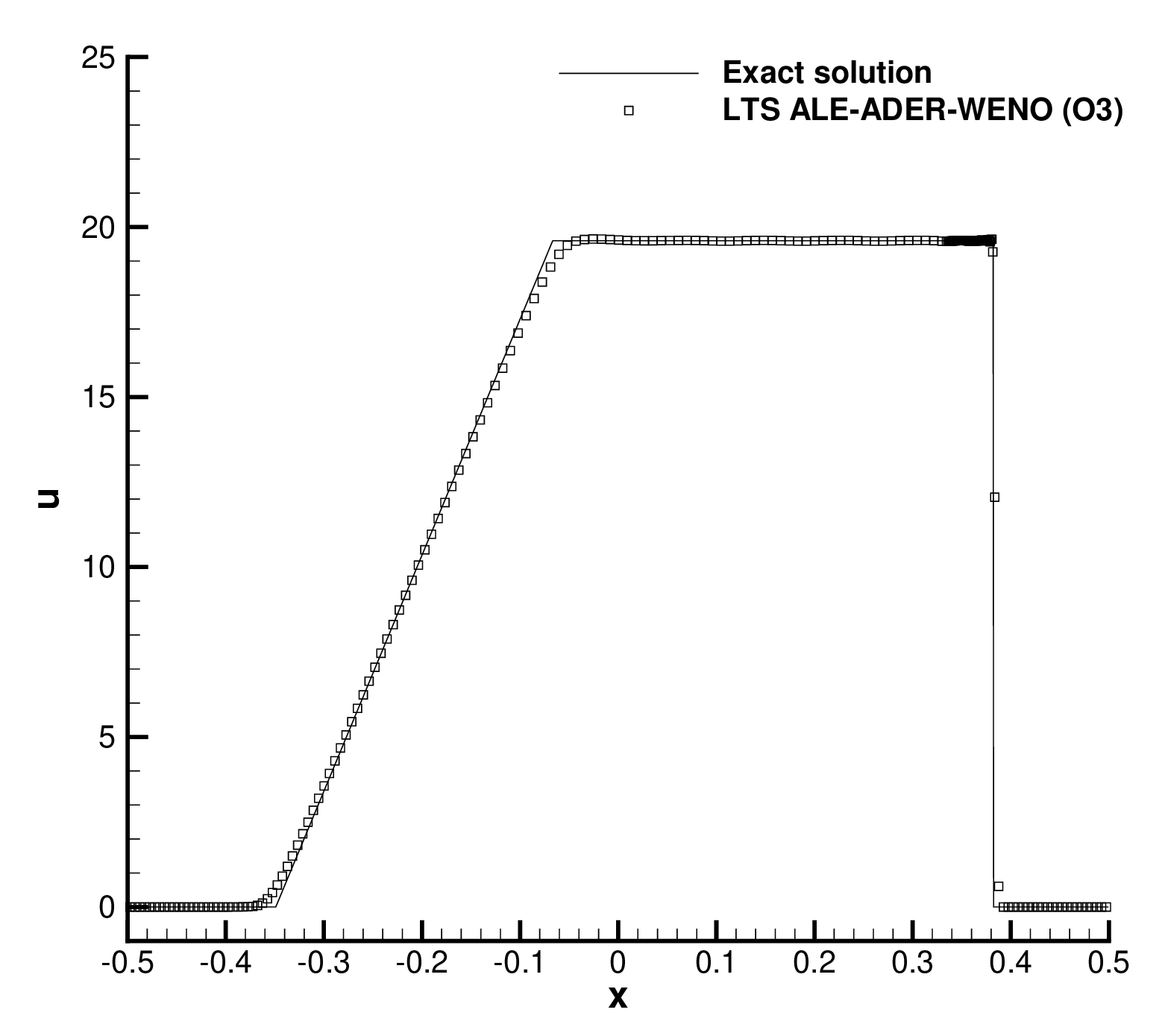}  \\
\multicolumn{2}{c}{\includegraphics[width=0.75\textwidth]{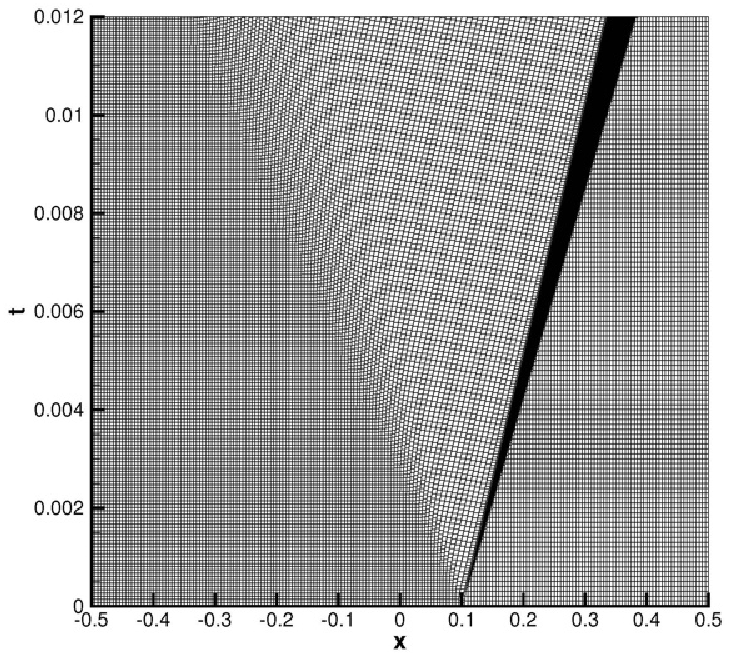}} 
\end{tabular}
\caption{Exact and numerical solution for shock tube problem RP3. Density (top left), velocity (top right) and resulting space-time mesh of the third order Lagrangian 
ADER-WENO finite volume scheme with conservative and time-accurate local time stepping (bottom).}
\label{fig.toro3}
\end{center}
\end{figure}

\begin{figure}
\begin{center}
\begin{tabular}{lr}
\includegraphics[width=0.45\textwidth]{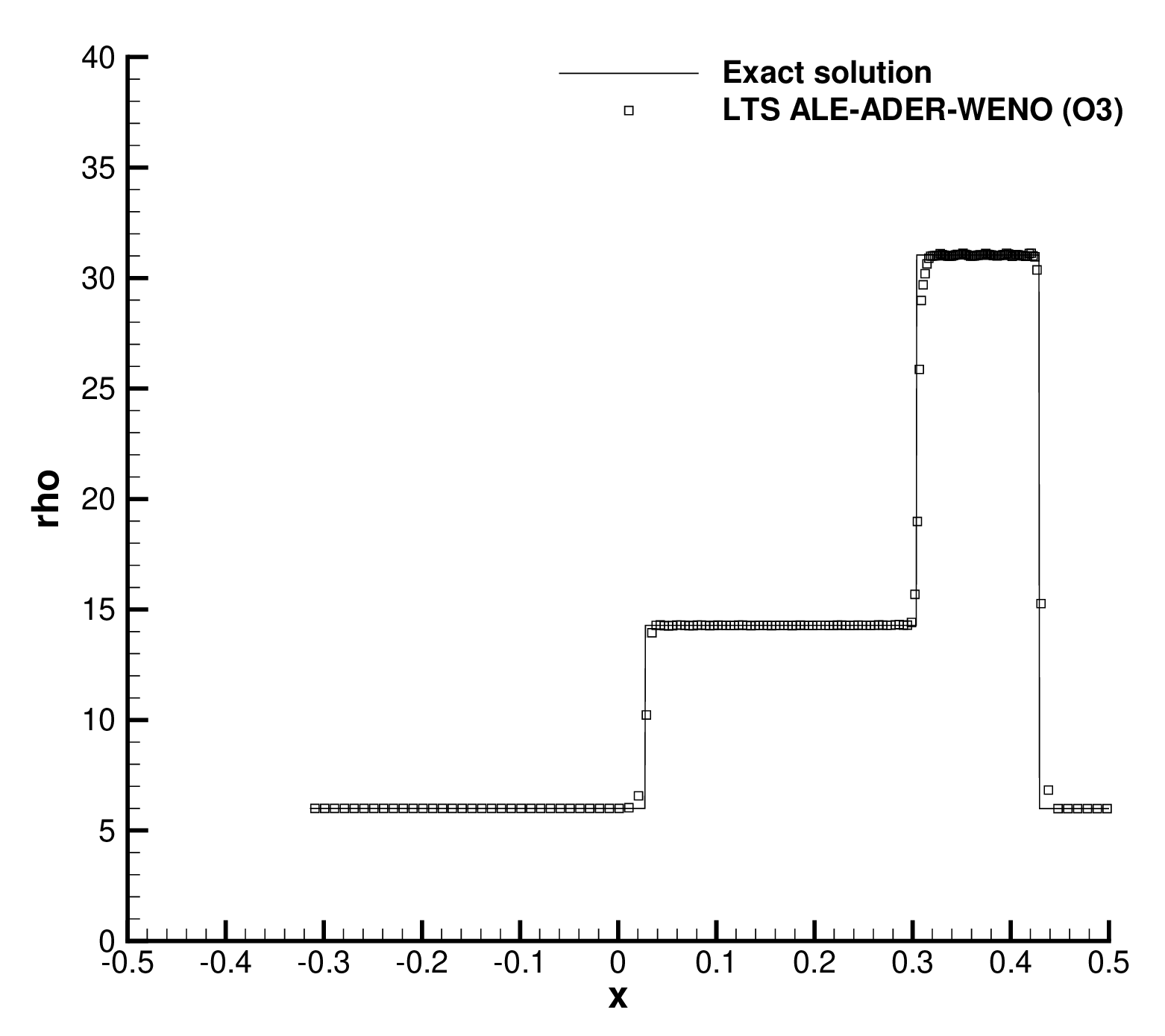} & 
\includegraphics[width=0.45\textwidth]{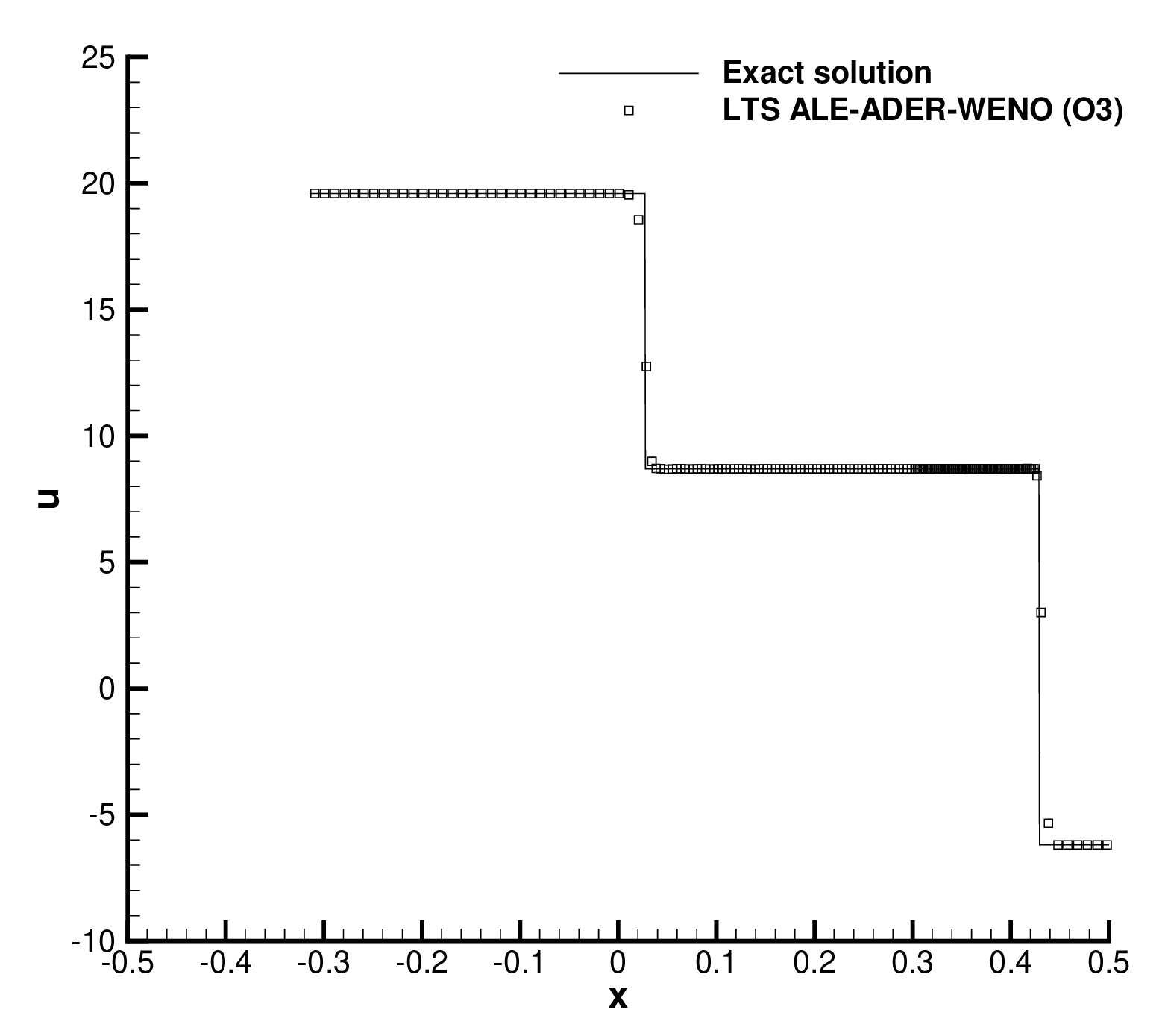}  \\
\multicolumn{2}{c}{\includegraphics[width=0.75\textwidth]{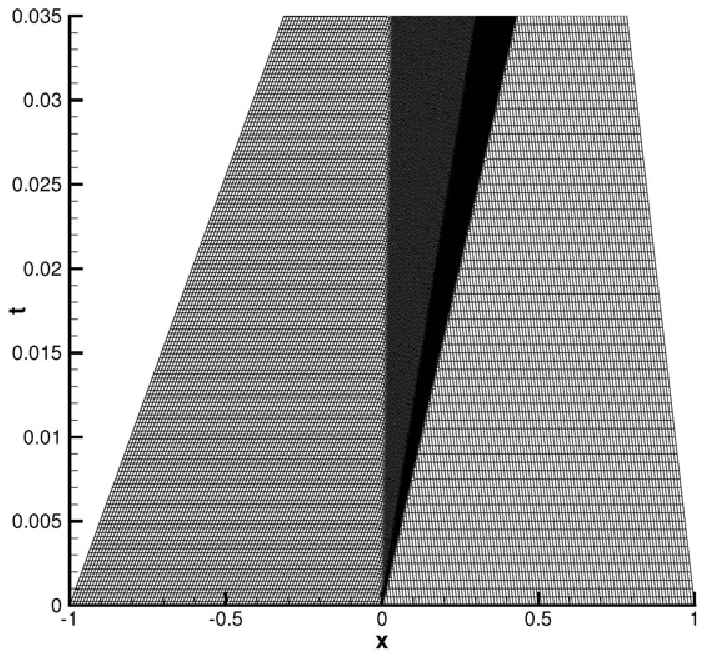}} 
\end{tabular}
\caption{Exact and numerical solution for shock tube problem RP4. Density (top left), velocity (top right) and resulting space-time mesh of the third order Lagrangian 
ADER-WENO finite volume scheme with conservative and time-accurate local time stepping (bottom).}
\label{fig.toro4}
\end{center}
\end{figure}

The initial conditions for the four Riemann problems RP1 - RP4 solved in the following are all of the type, 
\begin{equation}
  \Q(x,0) = \left\{ \begin{array}{ccc} \Q_L & \textnormal{ if } & x \leq x_d, \\ 
                                       \Q_R & \textnormal{ if } & x > x_d.         
                      \end{array}  \right. 
\end{equation}
The initial left and right states are listed in detail in Table \ref{tab.ic.euler}, together with the initial position of the discontinuity $x_d$ and the 
final simulation time $t_{\textnormal{end}}$. The initial computational  
domain is $\Omega(0)=[x_L;x_R]$, where $x_L$ and $x_R$ can be easily identified for each problem from the bottom of Figures \ref{fig.sod}-\ref{fig.toro4}; 
in all cases $\gamma=1.4$. RP1 and RP2 are the classical Sod and Lax shock tube problems, respectively, while RP3 and RP4 are taken from \cite{toro-book}.  
The exact solution has been computed using the exact Riemann solver described in great detail in \cite{toro-book}, where also a full FORTRAN listing is 
available. 

The computational results obtained for density $\rho$ and velocity $u$ with a third order ADER-WENO scheme (degree of the reconstruction polynomial $M=2$) 
are depicted together with the exact solution in the top row of Figures \ref{fig.sod}-\ref{fig.toro4}. The Osher-type flux \eqref{eqn.osher} is 
used in all cases. The initial spatial mesh consists of 200 initially equidistant control volumes for all test problems and the local CFL number 
is chosen as CFL=0.5. 
From the computational results one can observe that the contact wave is very well resolved and that the results are essentially non-oscillatory. 
Note that the scheme presented in this article is \textit{not} purely Lagrangian, but it is of the direct ALE type, hence a non-vanishing mass flux 
through the contact wave is possible, which leads to a small, but visible, smearing of the contact discontinuity (see Figures \ref{fig.sod}-\ref{fig.toro4}). 
The resulting space-time mesh for all test cases is shown in the bottom 
row of Figures \ref{fig.sod}-\ref{fig.toro4}. One can clearly identify the emerging wave structure of the Riemann problem, as well as the local time 
steps, which are given by the vertical extent of each single space-time control volume. Furthermore, the space-time meshes also show the hanging nodes
in time that arise within the LTS algorithm. For all cases RP1-RP4 it has been explicitly verified that the \textit{relative conservation error} 
was of the order of machine accuracy for all conserved quantities. The total absolute conservation error at the final simulation time is listed for each 
test problem and for all conserved quantities in Table \ref{tab.eu.euler}, together with the necessary total number of element updates when using a local 
time-stepping algorithm (LTS) or a traditional method based on global time stepping (GTS). 
The last column of Table \ref{tab.eu.euler} contains the savings factor, which clearly indicates that for the test problems considered here, a substantial 
savings in element updates can be achieved by using an LTS strategy. 
 
\begin{table}[!t]
 \caption{Initial states left and right for the density $\rho$, velocity $u$ 
 and the pressure $p$ for the compressible Euler equations. The final output times,  ($t_{\textnormal{end}}$) 
 and the initial position of the discontinuity ($x_d$) are also given. } 
\begin{center} 
 \begin{tabular}{ccccccccc}
 \hline
 Case & $\rho_L$ & $u_L$ & $p_L$ & $\rho_R$ & $u_R$ & $p_R$ & $t_{\textnormal{end}}$ & $x_d$        \\ 
 \hline
 RP1    &  1.0      &  0.0       & 1.0     & 0.125      &  0.0        & 0.1      & 0.4    &  0.0    \\
 RP2    &  0.445    &  0.698     & 3.528   & 0.5        &  0.0        & 0.571    & 0.1    &  0.0    \\
 RP3    &  1.0      &  0.0       & 1000    & 1.0        &  0.0        & 0.01     & 0.012  &  0.1    \\
 RP4    &  5.99924  &  19.5975   & 460.894 & 5.99242    & -6.19633    & 46.095   & 0.035  & -0.2    \\ 
 \hline
 \end{tabular}
\end{center} 
 \label{tab.ic.euler}   
\end{table}

\begin{table}[!b]
 \caption{Euler equations: Total absolute conservation error for each test case and comparison of the computational efficiency between GTS and LTS algorithm using the total 
          number of element updates. } 
\begin{center} 
 \begin{tabular}{ccccccc}
 \hline
	    & \multicolumn{3}{c}{Conservation error} & \multicolumn{3}{c}{Number of element updates} \\ 
 Case & $\rho$ & $\rho u$ & $\rho E$ & GTS & LTS & GTS/LTS        \\ 
 \hline
 RP1    &  2.08944E-13  &  4.60021E-13  & 9.86766E-13   & 73600      & 21840  & 3.37  \\
 RP2    &  2.24598E-13  &  1.69381E-12  & 7.67386E-12   & 65000      & 20579  & 3.16  \\
 RP3    &  1.13805E-11  &  9.82959E-10  & 5.47216E-08   & 167000     & 34222  & 4.88  \\
 RP4    &  1.23741E-11  &  3.40464E-10  & 5.94946E-09   & 294800     & 50481  & 5.84  \\ 
 \hline
 \end{tabular}
\end{center} 
 \label{tab.eu.euler}   
\end{table}

\subsection{Ideal classical MHD equations} 
\label{sec.mhd}

In this section the equations of ideal classical magnetohydrodynamics (MHD) are solved with the proposed Arbitrary-Lagrangian-Eulerian ADER-WENO
finite volume scheme, together with the new local time stepping (LTS) feature. The augmented PDE system including the hyperbolic divergence-correction 
term proposed by Dedner et al. \cite{dedneretal} reads 
\begin{equation}
\label{eqn.MHD} 
  \frac{\partial}{\partial t} \left( \begin{array}{c} 
  \rho \\   \rho \vec v \\   \rho E \\ \vec B \\ \psi \end{array} \right)   
  + \nabla \cdot \left( \begin{array}{c} 
  \rho \vec v  \\ 
  \rho \vec v \vec v + p_t \Id - \frac{1}{4 \pi} \vec B \vec B \\ 
  \vec v (\rho E + p_t ) - \frac{1}{4 \pi} \vec B ( \vec v \cdot \vec B ) \\ 
  \vec v \vec B - \vec B \vec v + \psi \Id \\
  c_h^2 \vec B 
  \end{array} \right) = 0, 
\end{equation} 
with 
\begin{equation}
  p = (\gamma - 1) (\rho E - \halb \rho \vec v^2 - \frac{1}{8 \pi} \vec B ^2), \qquad p_t = p + \frac{1}{8 \pi} \vec B ^2.
\end{equation} 
In Eqn. \ref{eqn.MHD}, $\rho$ is the density of the gas, $\vec v = (u,v,w)$ is the velocity vector, $\vec B = (B_x,B_y,B_z)$ is the 
vector of the magnetic field, $p$ is the gas pressure, $p_t$ is the sum of the gas and the magnetic pressure, $\rho E$ is the 
total energy density and $\gamma$ is the ratio of specific heats. $\Id$ is the unit matrix and the notation $\vec x_1 \vec x_2$ denotes
the dyadic product of two vectors $\vec x_1$ and $\vec x_2$. The scalar $\psi$ is used for divergence cleaning, see \cite{dedneretal},  
to satisfy the constraint $\nabla \cdot \vec B = 0$. In one space dimension, this constraint simply reduces to $\partial B_x / \partial x = 0$. 

In the following one-dimensional test problems either the classical Rusanov flux \eqref{eqn.rusanov} is used, or the new Osher-type flux 
given by Eqn. \eqref{eqn.osher}, which has been originally proposed by Dumbser and Toro for the Eulerian case in \cite{OsherUniversal,OsherNC} 
and has subsequently been extended to Lagrangian schemes in \cite{Lagrange1D}. 

\begin{figure}
\begin{center}
\begin{tabular}{lr}
\includegraphics[width=0.45\textwidth]{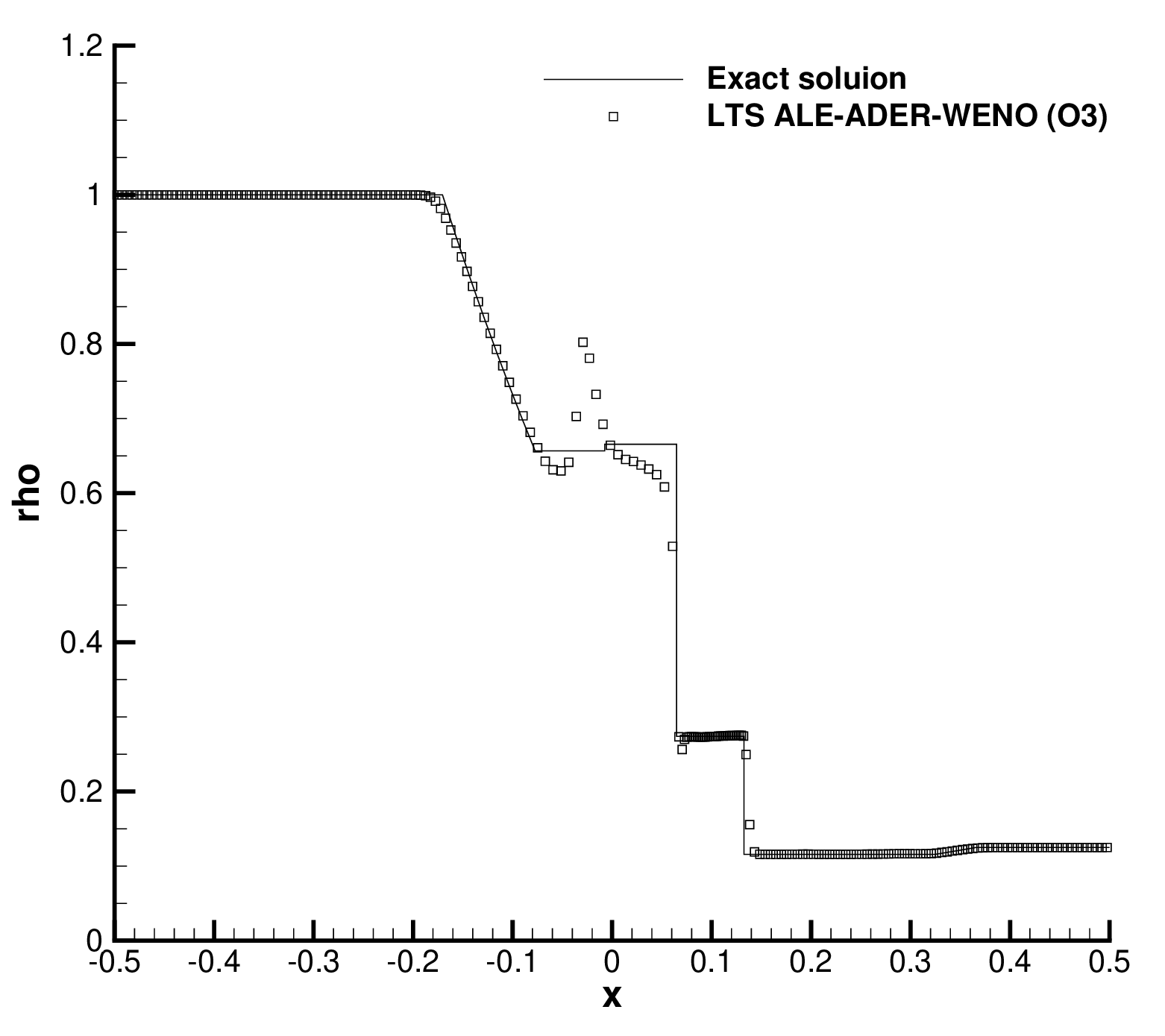} & 
\includegraphics[width=0.45\textwidth]{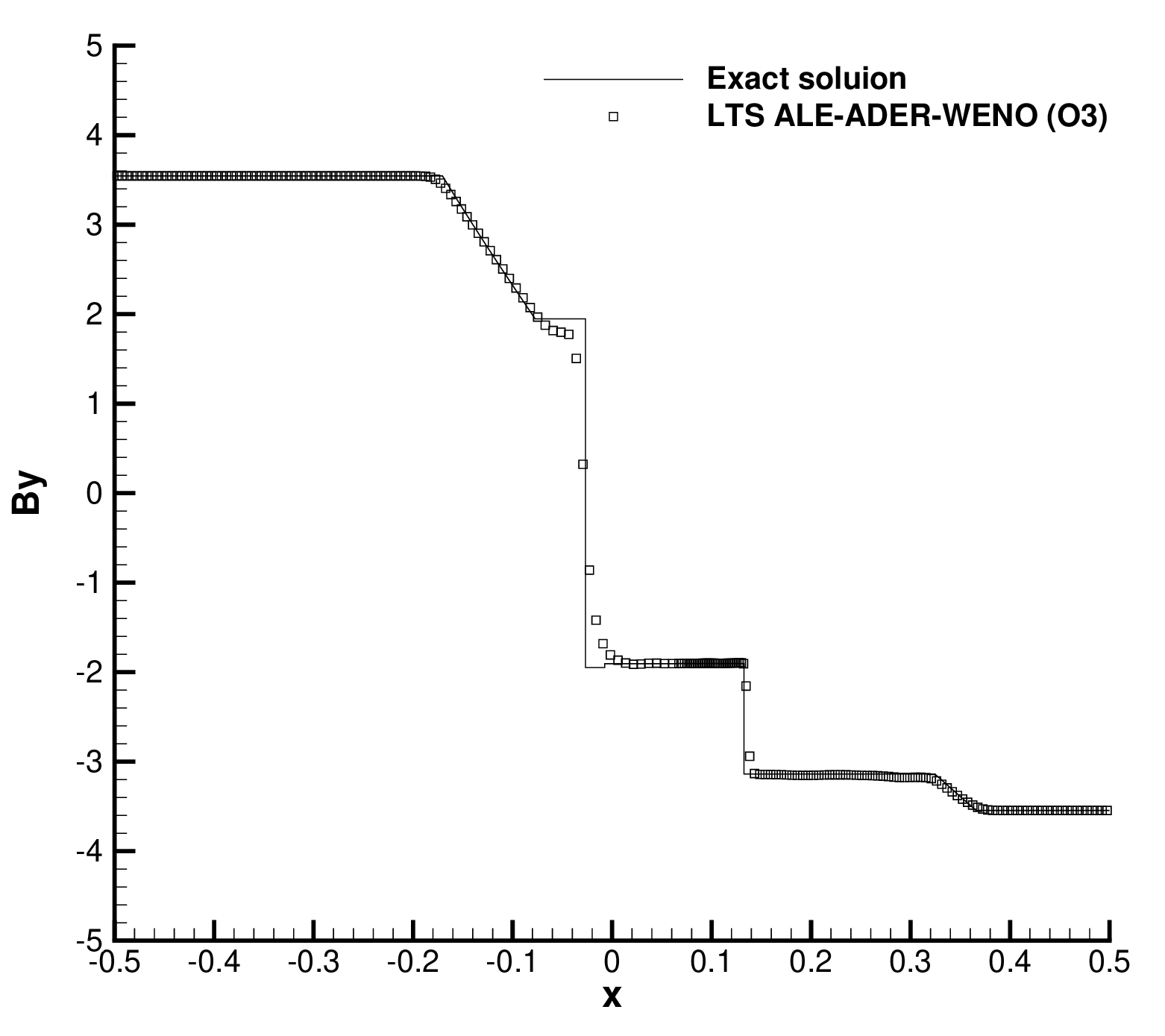}  \\
\multicolumn{2}{c}{\includegraphics[width=0.75\textwidth]{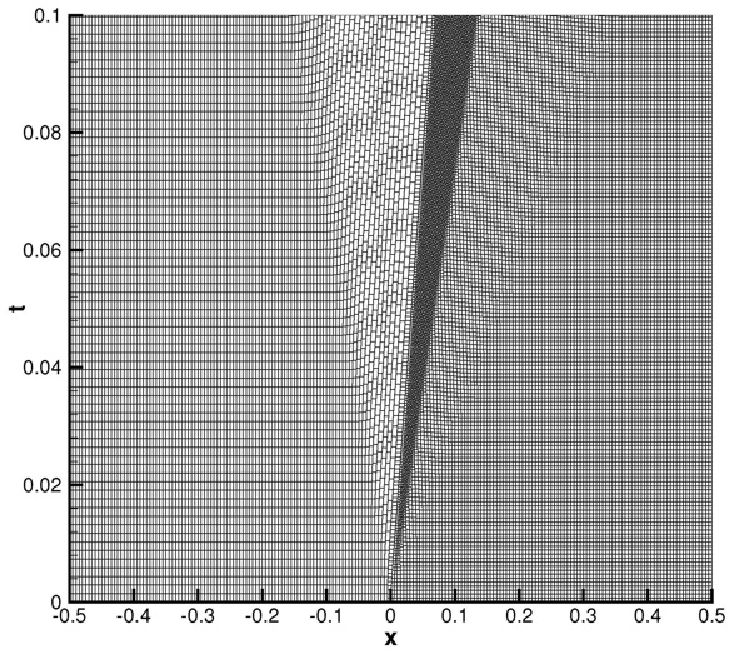}} 
\end{tabular}
\caption{Exact and numerical solution for the MHD shock tube problem RP1. Density (top left), magnetic field component $B_y$ (top right) and resulting space-time 
mesh of the third order Lagrangian ADER-WENO finite volume scheme with conservative and time-accurate local time stepping (bottom).}
\label{fig.mhd1}
\end{center}
\end{figure}

\begin{figure}
\begin{center}
\begin{tabular}{lr}
\includegraphics[width=0.45\textwidth]{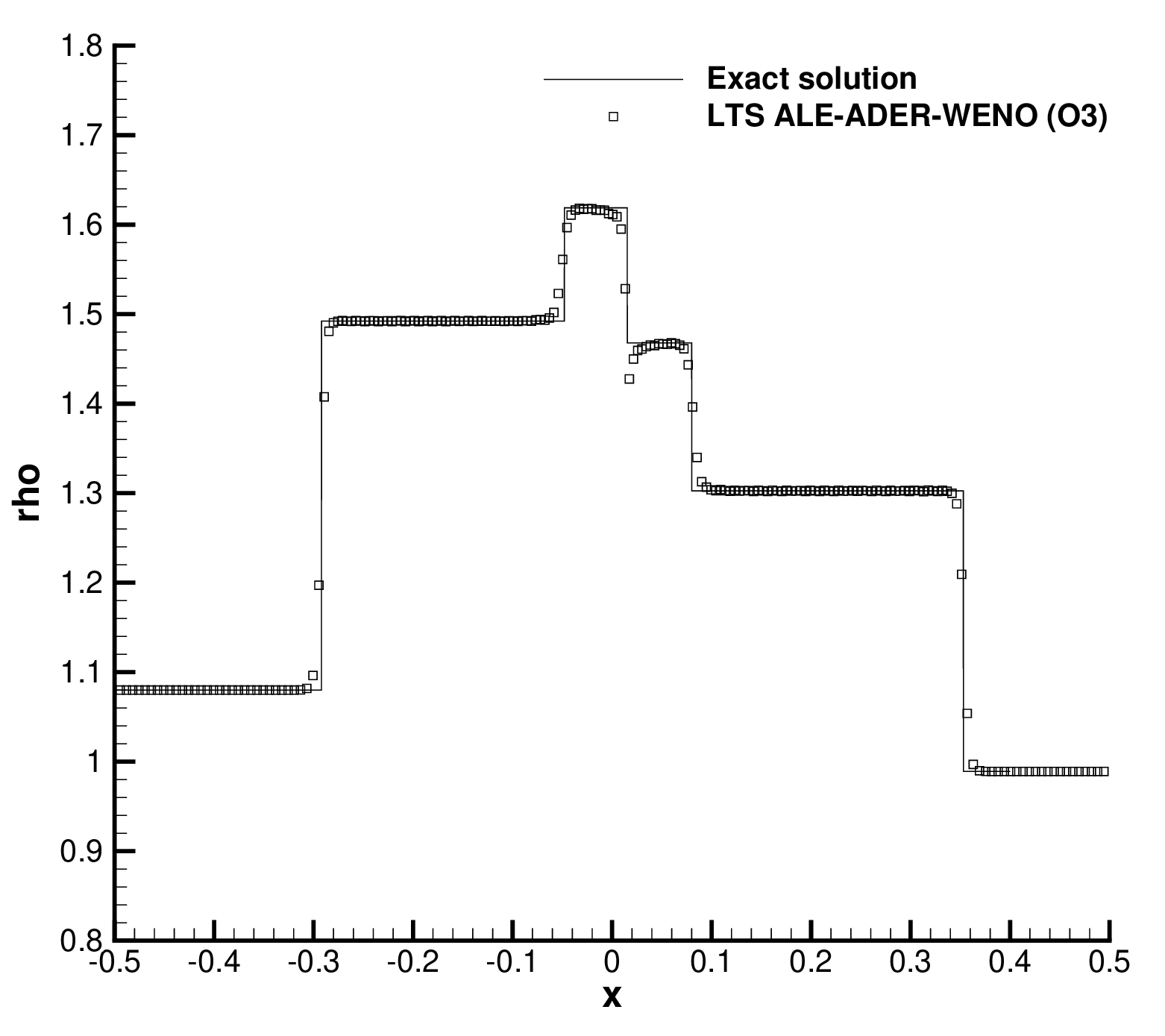} & 
\includegraphics[width=0.45\textwidth]{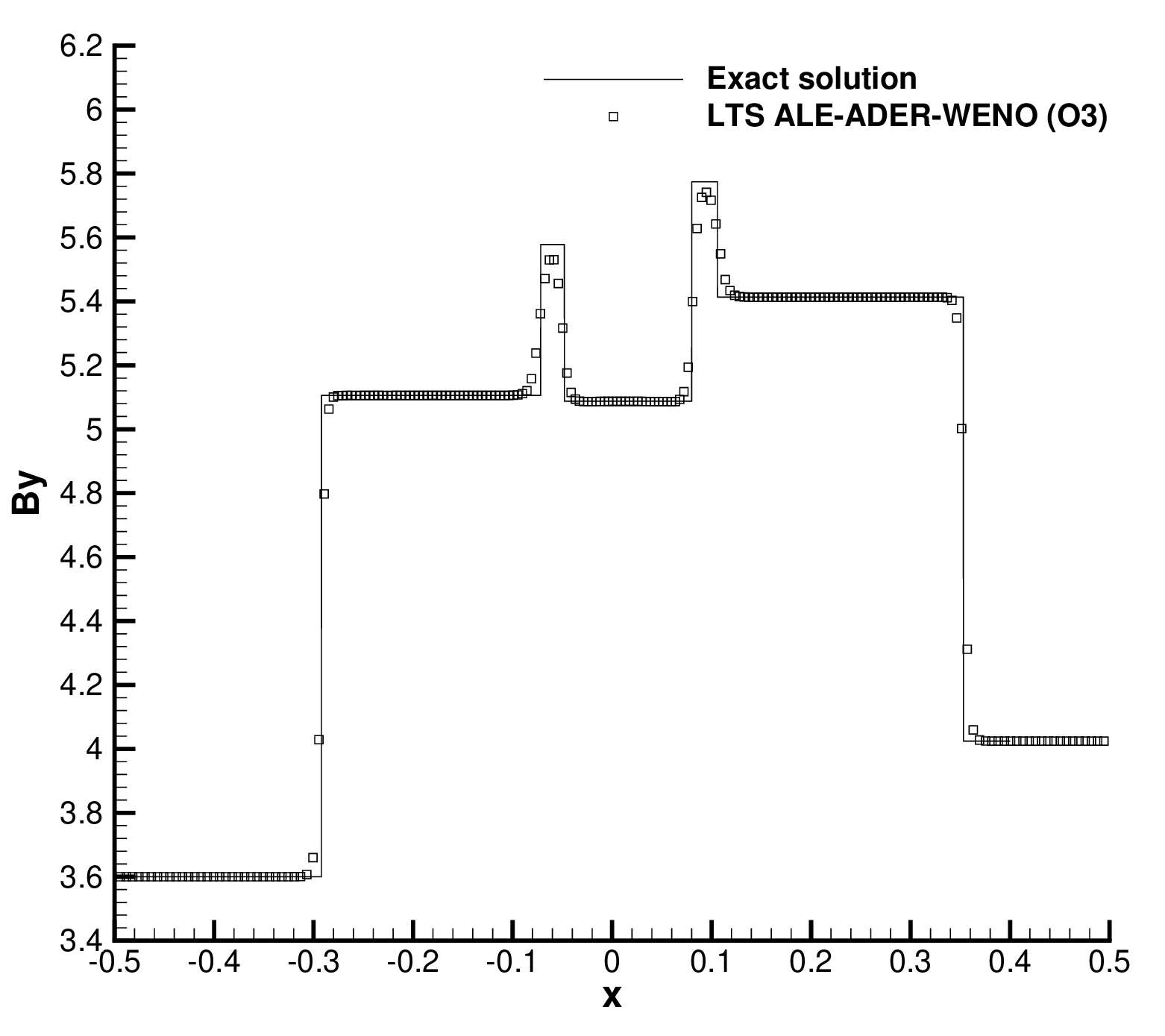}  \\
\multicolumn{2}{c}{\includegraphics[width=0.75\textwidth]{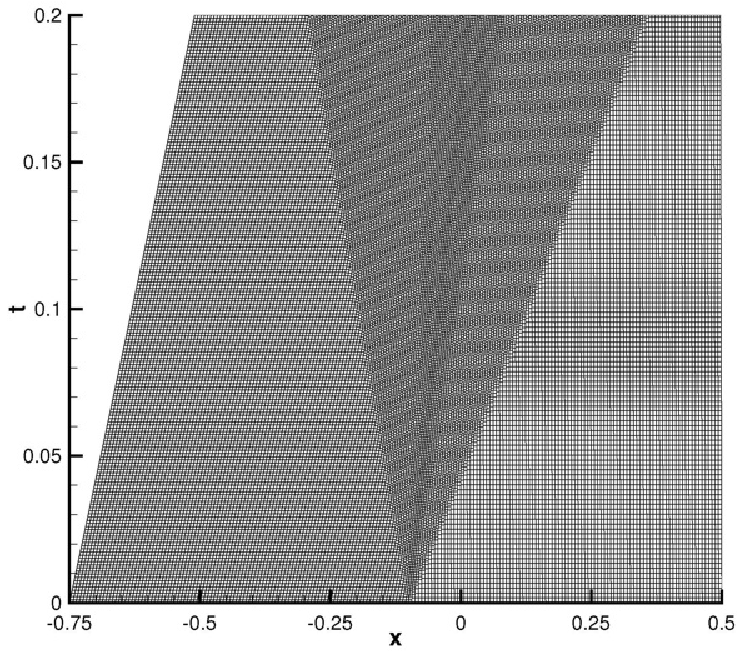}} 
\end{tabular}
\caption{Exact and numerical solution for the MHD shock tube problem RP2. Density (top left), magnetic field component $B_y$ (top right) and resulting space-time 
mesh of the third order Lagrangian ADER-WENO finite volume scheme with conservative and time-accurate local time stepping (bottom).}
\label{fig.mhd2}
\end{center}
\end{figure}

\begin{figure}
\begin{center}
\begin{tabular}{lr}
\includegraphics[width=0.45\textwidth]{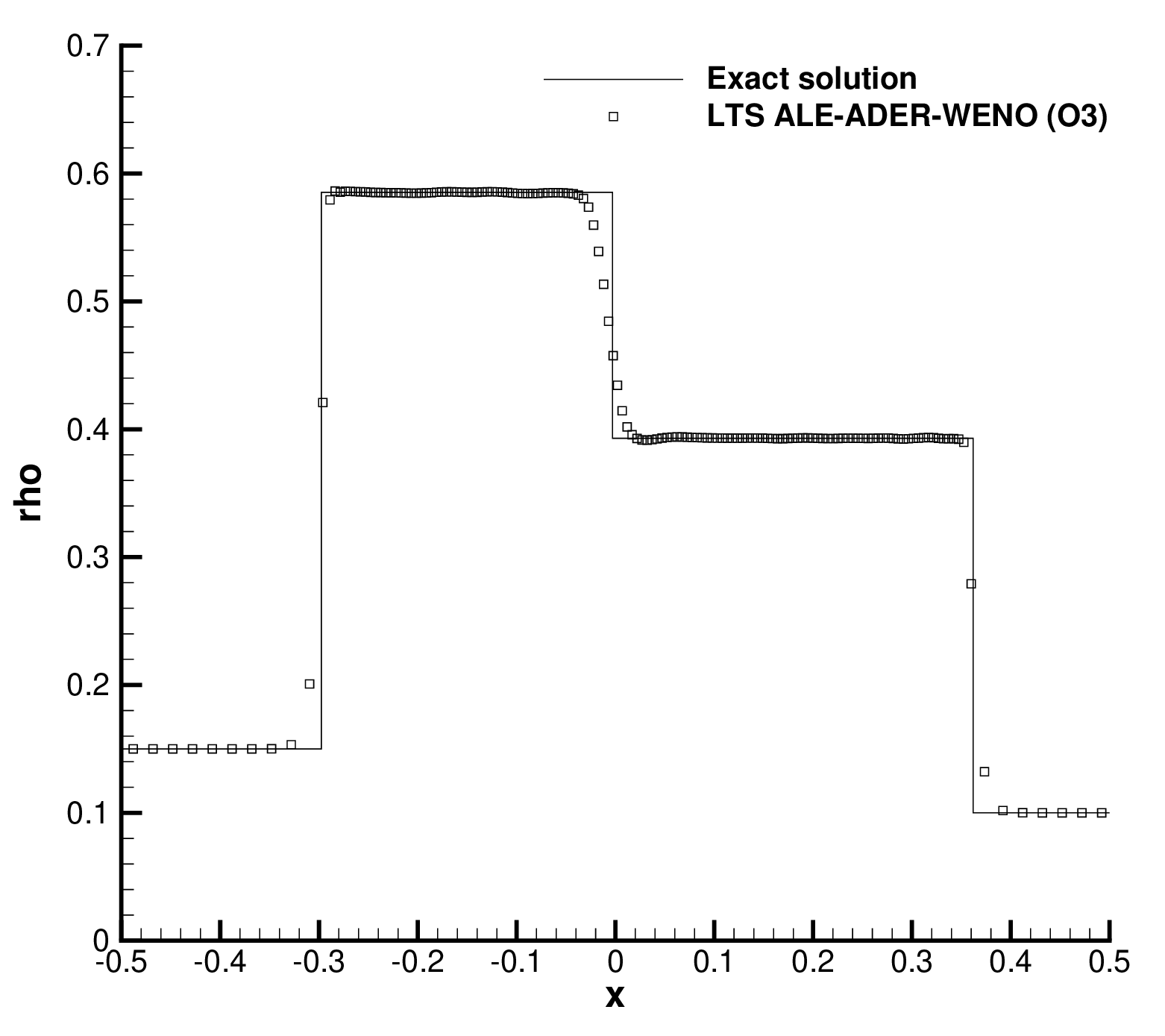} & 
\includegraphics[width=0.45\textwidth]{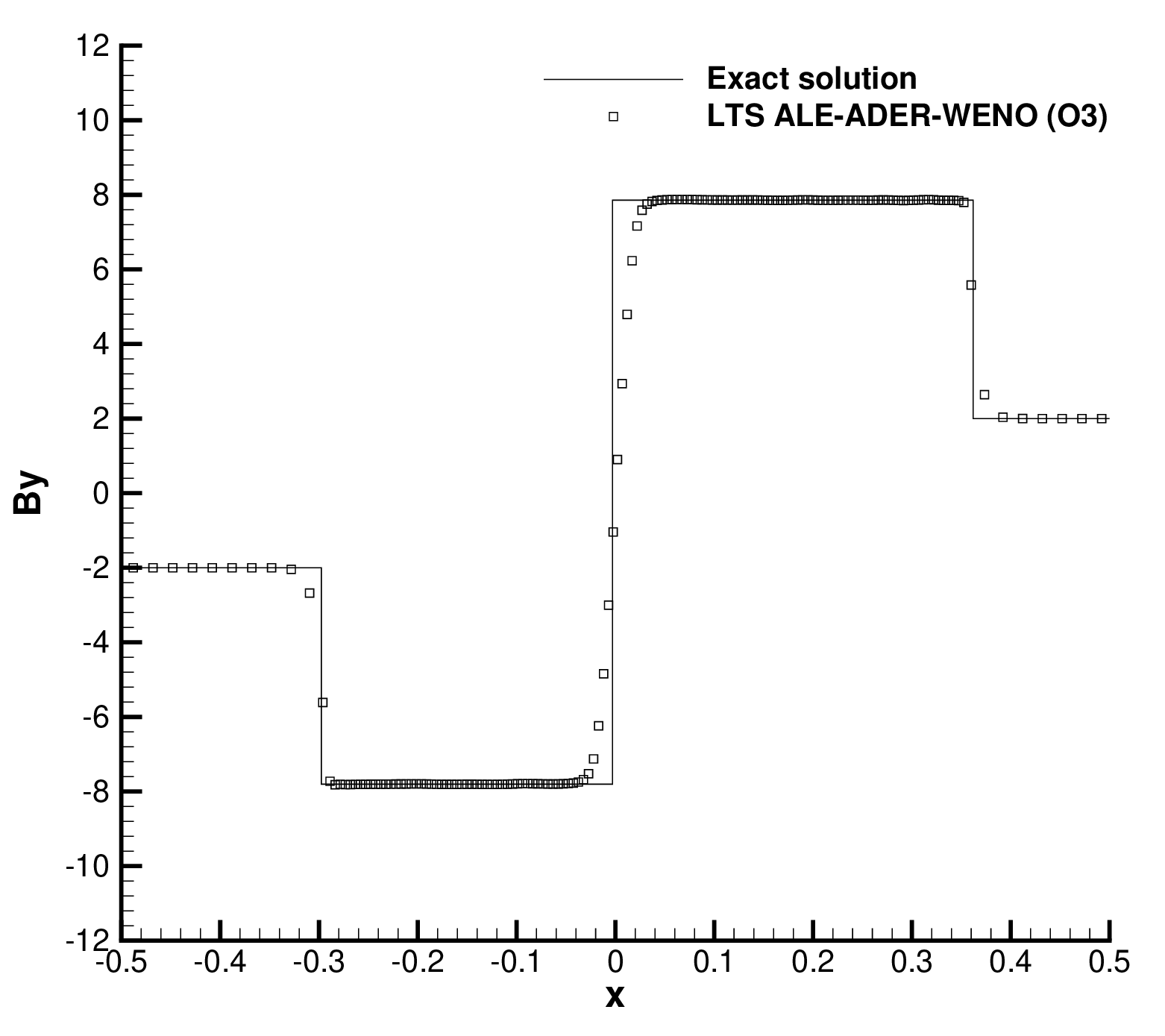}  \\
\multicolumn{2}{c}{\includegraphics[width=0.75\textwidth]{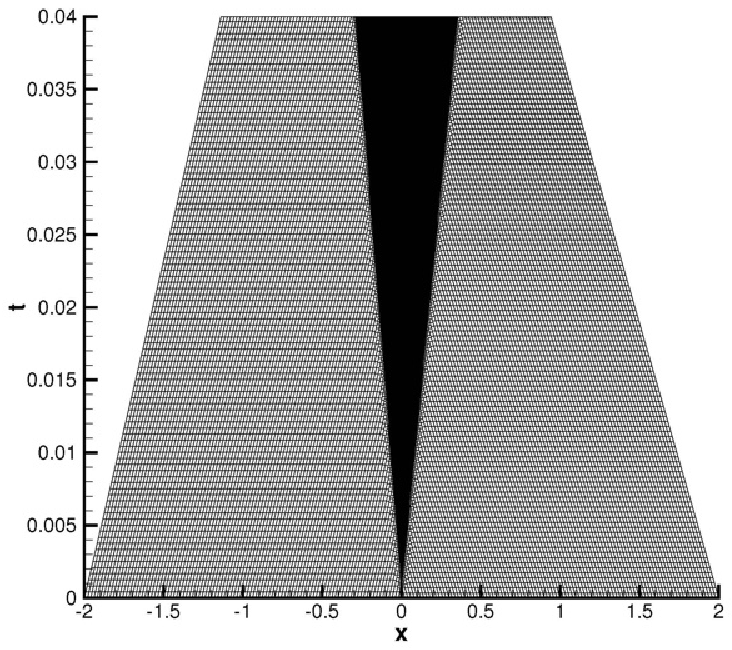}} 
\end{tabular}
\caption{Exact and numerical solution for the MHD shock tube problem RP3. Density (top left), magnetic field component $B_y$ (top right) and resulting space-time 
mesh of the third order Lagrangian ADER-WENO finite volume scheme with conservative and time-accurate local time stepping (bottom).}
\label{fig.mhd3}
\end{center}
\end{figure}

\begin{figure}
\begin{center}
\begin{tabular}{lr}
\includegraphics[width=0.45\textwidth]{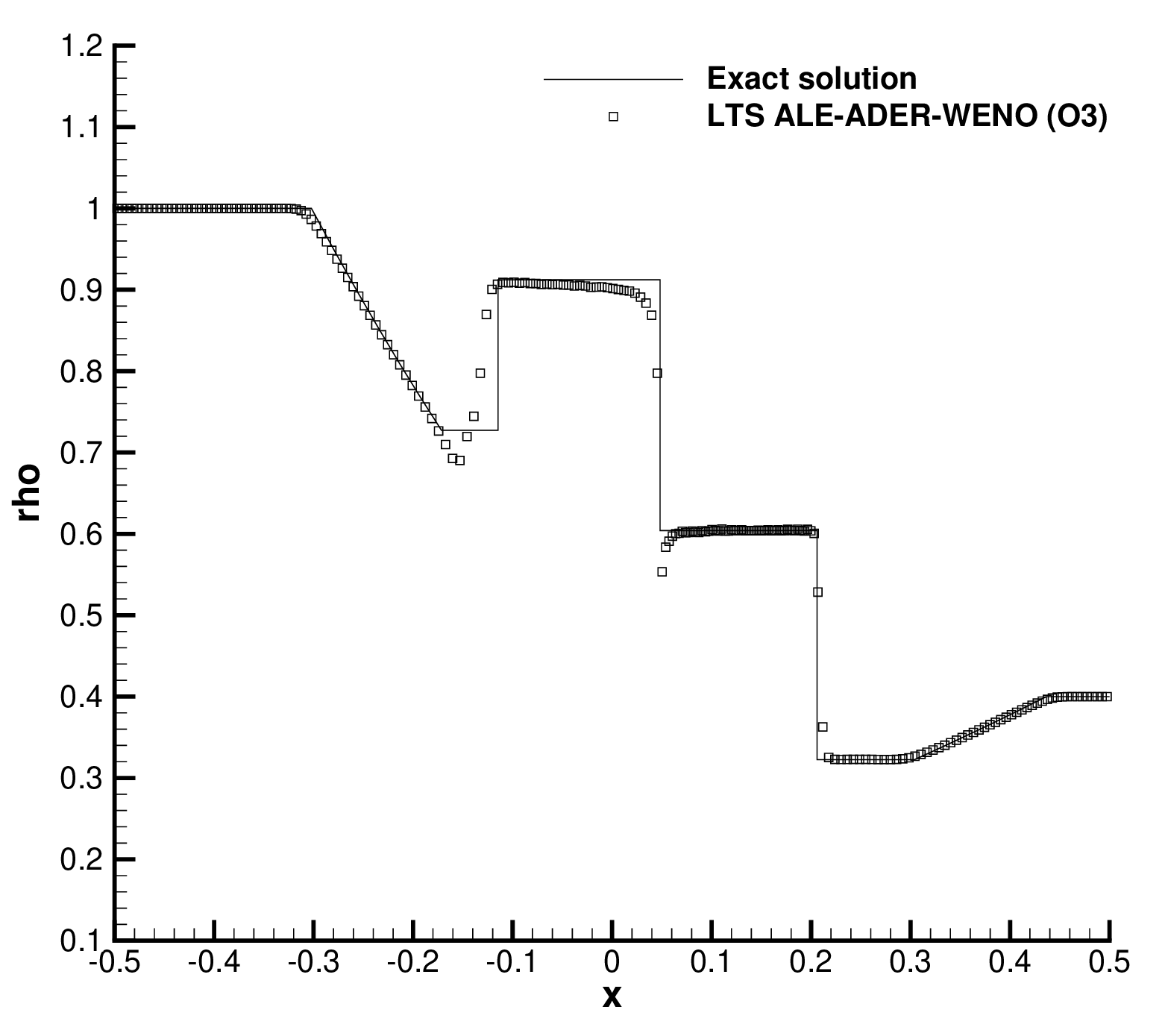} & 
\includegraphics[width=0.45\textwidth]{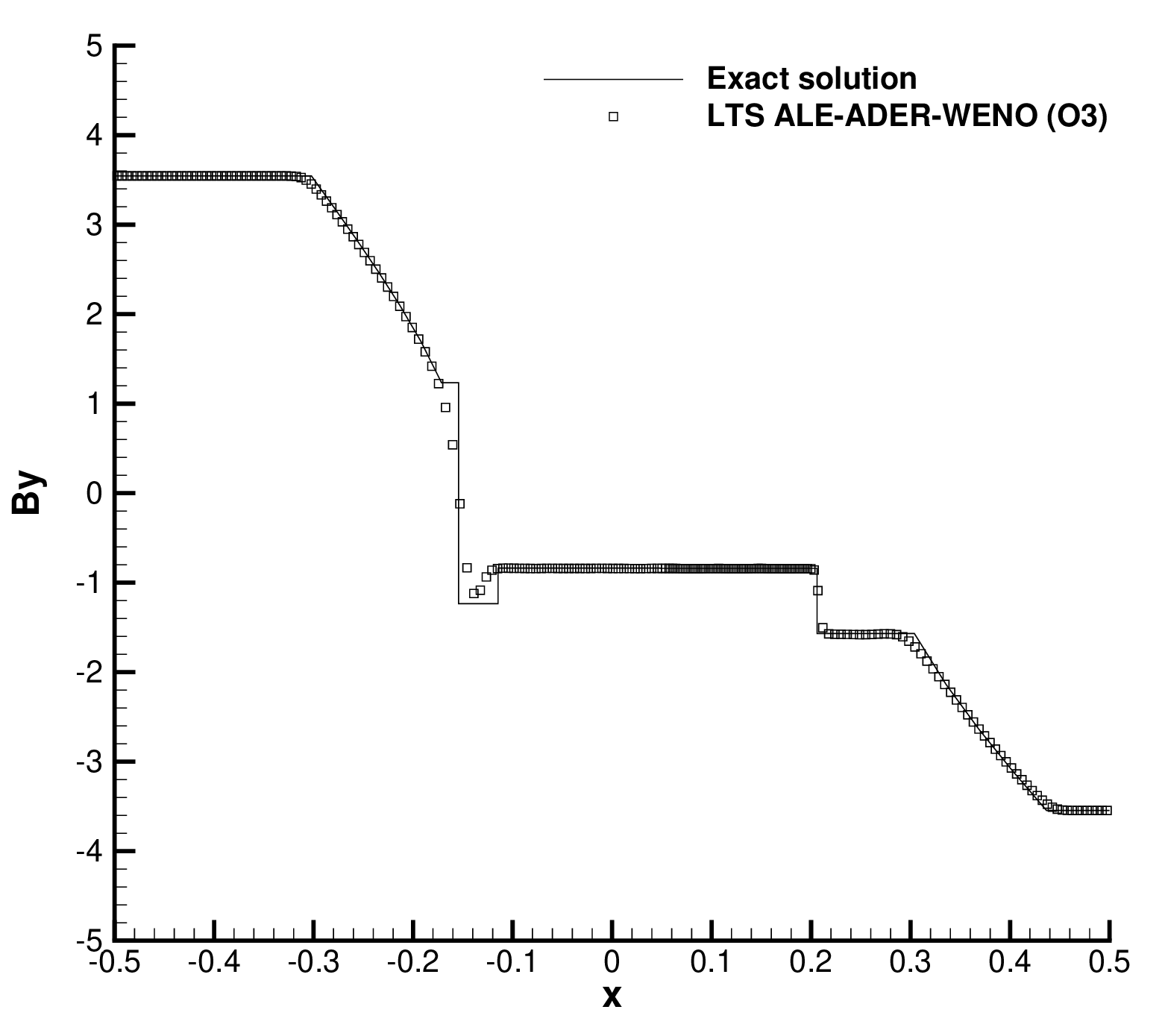}  \\
\multicolumn{2}{c}{\includegraphics[width=0.75\textwidth]{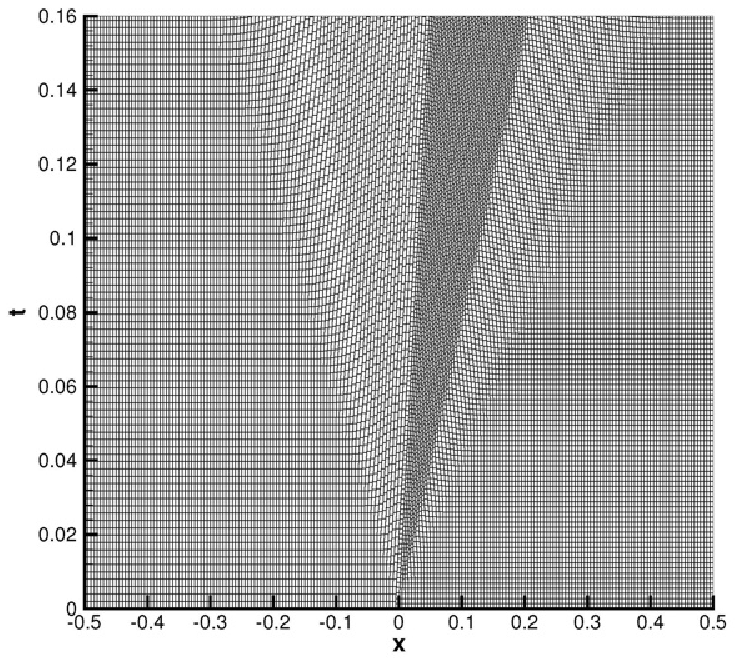}} 
\end{tabular}
\caption{Exact and numerical solution for the MHD shock tube problem RP4. Density (top left), magnetic field component $B_y$ (top right) and resulting space-time 
mesh of the third order Lagrangian ADER-WENO finite volume scheme with conservative and time-accurate local time stepping (bottom).}
\label{fig.mhd4}
\end{center}
\end{figure}

\begin{figure}
\begin{center}
\begin{tabular}{lr}
\includegraphics[width=0.45\textwidth]{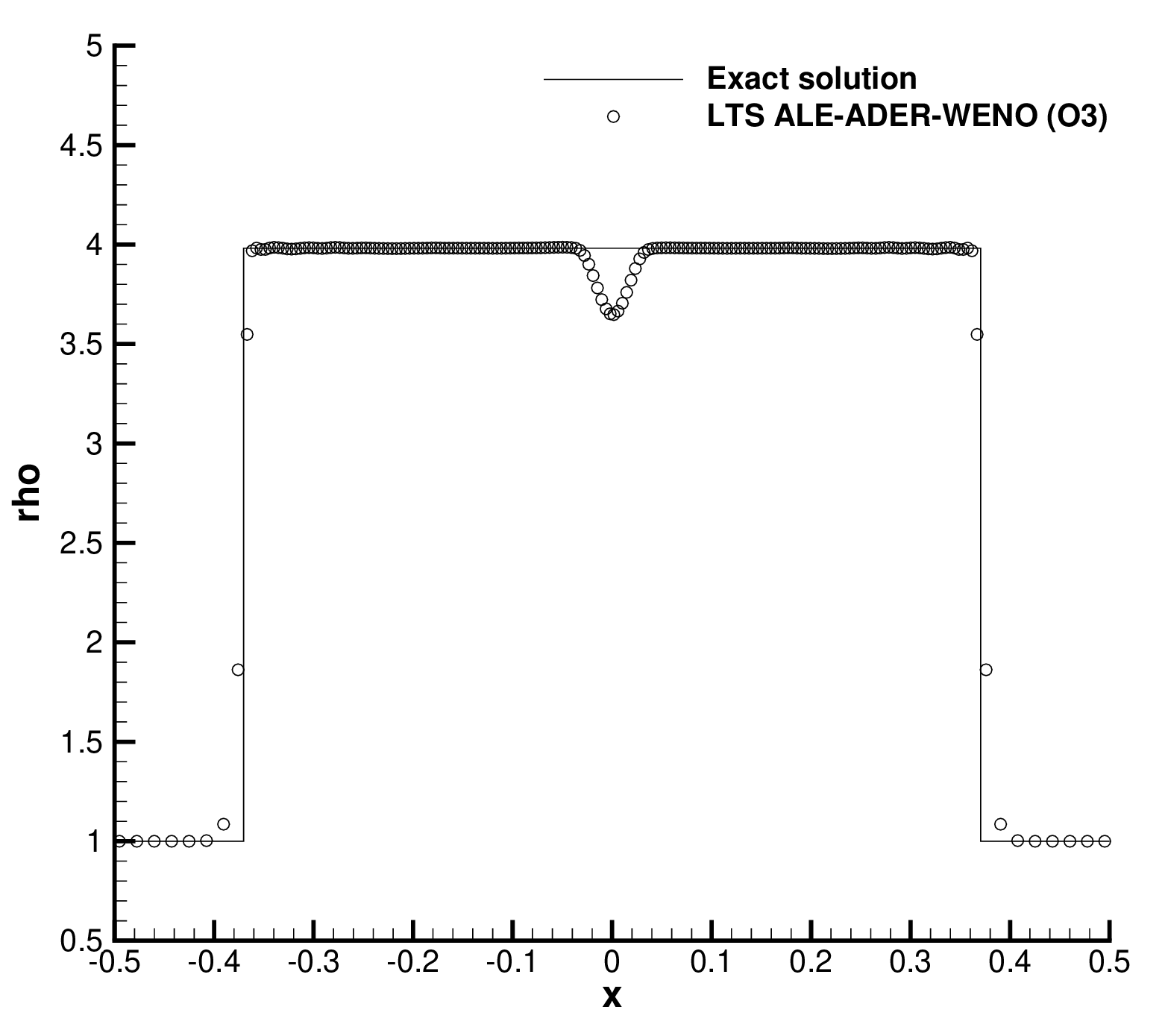} & 
\includegraphics[width=0.45\textwidth]{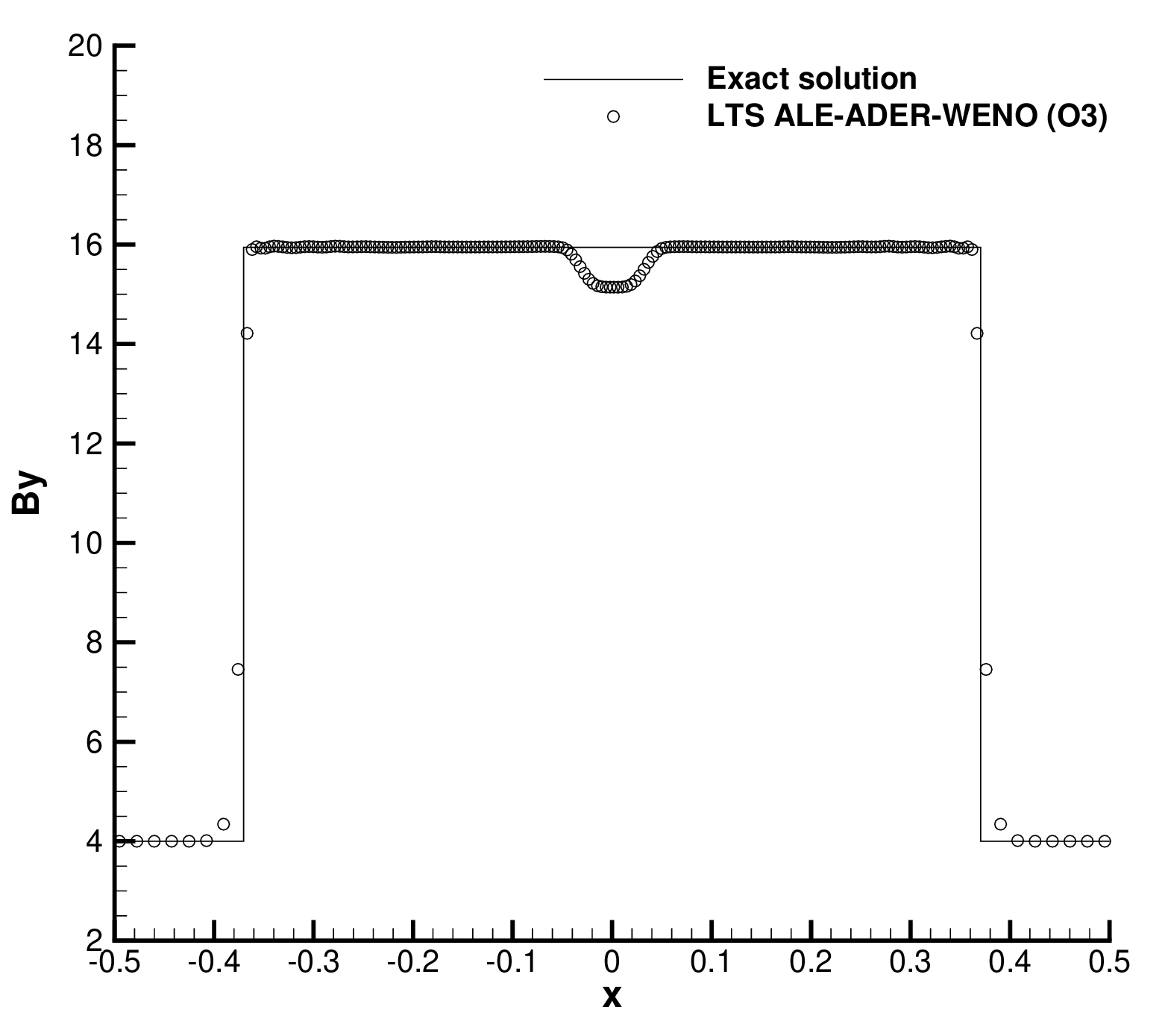}  \\
\multicolumn{2}{c}{\includegraphics[width=0.75\textwidth]{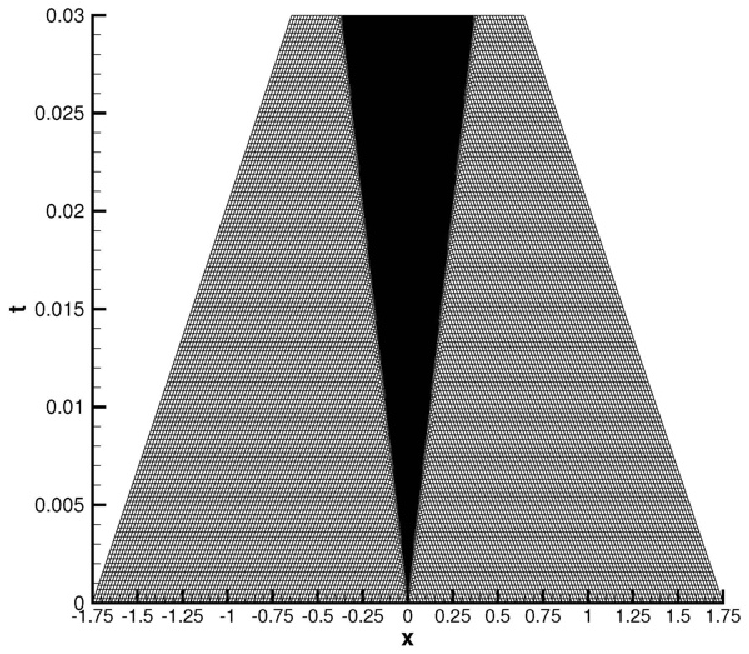}} 
\end{tabular}
\caption{Exact and numerical solution for the MHD shock tube problem RP5. Density (top left), magnetic field component $B_y$ (top right) and resulting space-time 
mesh of the third order Lagrangian ADER-WENO finite volume scheme with conservative and time-accurate local time stepping (bottom).}
\label{fig.mhd5}
\end{center}
\end{figure}

\begin{figure}
\begin{center}
\begin{tabular}{lr}
\includegraphics[width=0.45\textwidth]{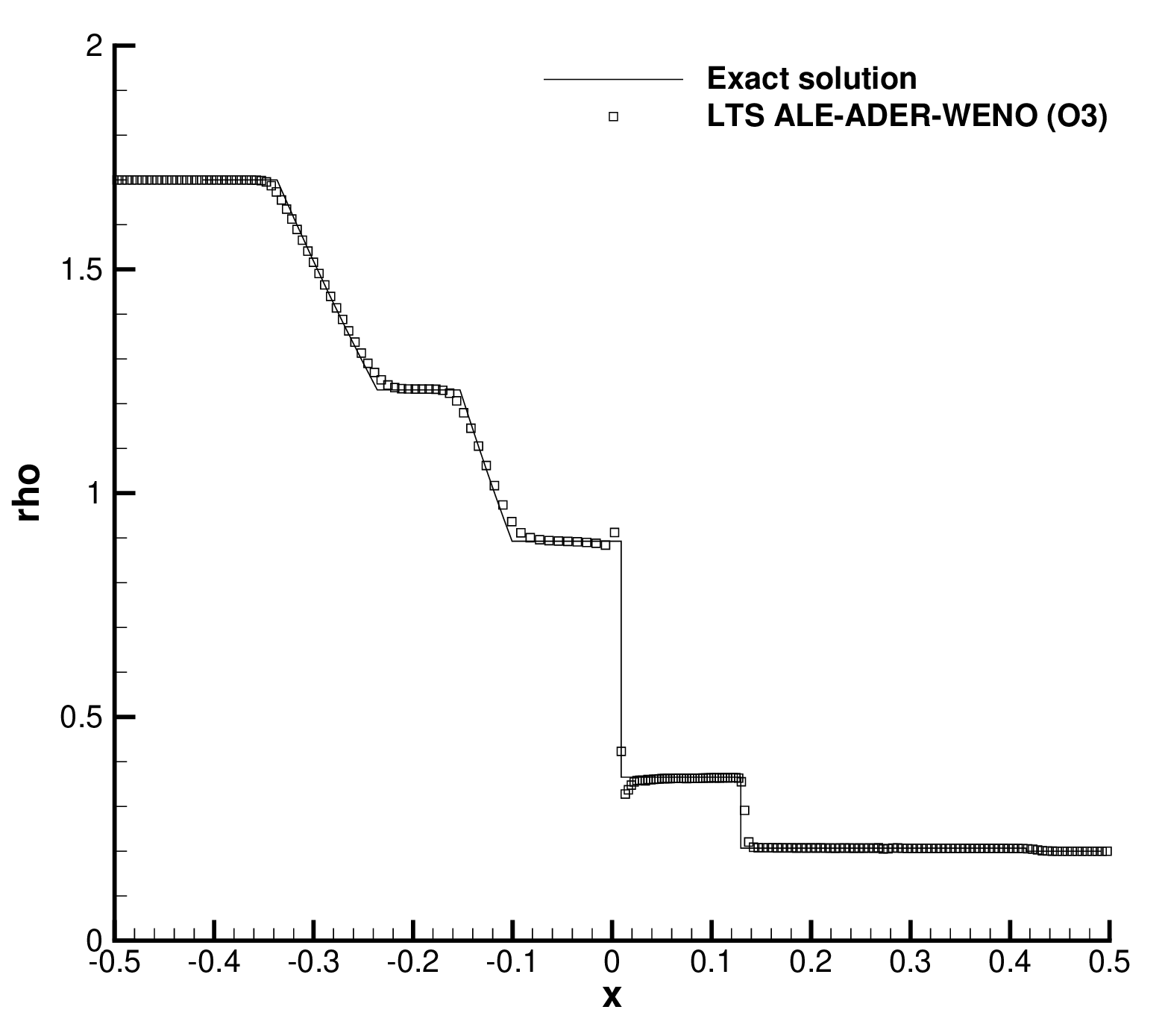} & 
\includegraphics[width=0.45\textwidth]{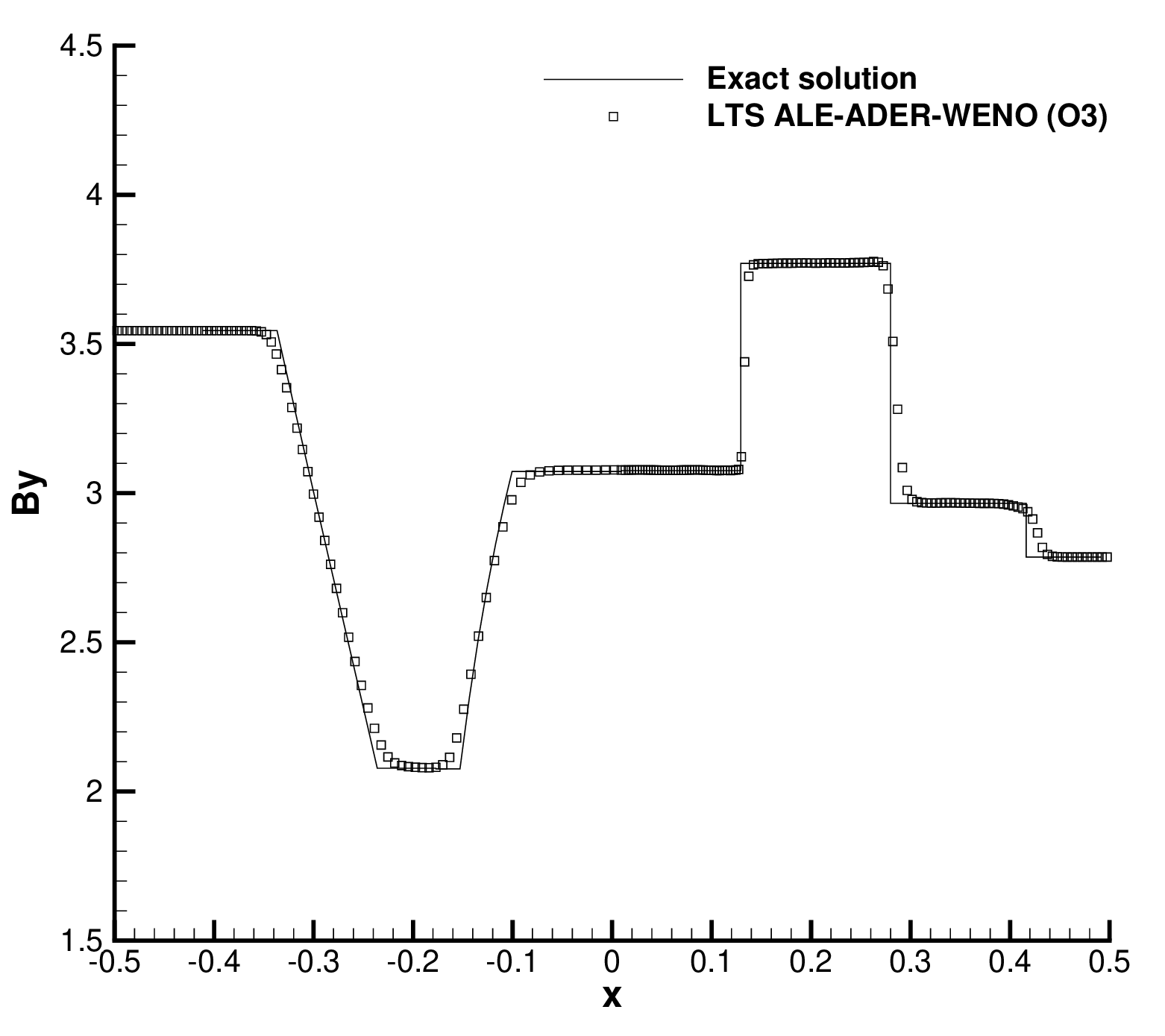}  \\
\multicolumn{2}{c}{\includegraphics[width=0.75\textwidth]{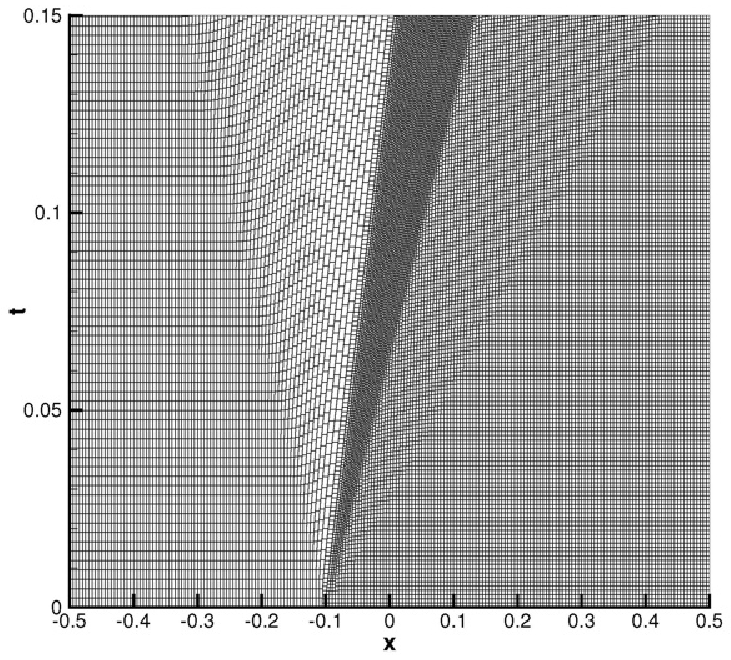}} 
\end{tabular}
\caption{Exact and numerical solution for the MHD shock tube problem RP6. Density (top left), magnetic field component $B_y$ (top right) and resulting space-time 
mesh of the third order Lagrangian ADER-WENO finite volume scheme with conservative and time-accurate local time stepping (bottom).}
\label{fig.mhd6}
\end{center}
\end{figure}

\paragraph{Riemann problems} 
The initial conditions for the shock tube problems are listed in Table \ref{tab.ic.mhd} and the ratio of specific heats is set to 
$\gamma=\frac{5}{3}$ for all cases. All shock-tube problems are solved on a mesh of 200 initially equidistant cells using the third 
order version of the ALE ADER-WENO scheme with LTS. 
The computational results are depicted in Figs. \ref{fig.mhd1} to \ref{fig.mhd6}, together with the exact solution and the resulting  
space-time meshes. The exact Riemann solver for MHD has kindly been provided by S.A.E.G. Falle \cite{fallemhd}. For an alternative exact 
Riemann solver of the MHD equations, see also \cite{torrilhon}.  
In all cases, the numerical results agree well with the exact solution, the contact wave is well resolved and most of the details of the 
wave structure emerging from the Riemann problem are properly resolved. As can be seen from Figures \ref{fig.mhd1}-\ref{fig.mhd6} the 
space-time mesh directly reflects the wave pattern of the Riemann problem, as expected. One can furthermore observe that in those regions where 
the mesh is compressed by the presence of a shock wave, the time step is locally decreased, while inside the rarefaction fans, where the mesh is 
expanding, the time step is locally increased. To assess the efficiency of the proposed local time stepping algorithm, the total number of element
updates is reported for all test problems in Table \ref{tab.eu.mhd}. Furthermore, the total conservation error at the final time is reported for 
each test case in Table \ref{tab.eu.mhd}, showing that also with the use of an LTS strategy a conservative scheme can be designed in the framework
of high order Lagrangian finite volume schemes.

\begin{table}[!b]
 \caption{Initial states left and right for the density $\rho$, velocity vector $\vec v=(u,v,w)$,  
 the pressure $p$ and the magnetic field vector $\vec B = (B_x,B_y,B_z)$ for the ideal classical MHD equations. The final output times,   ($t_{\textnormal{end}}$) 
 and the initial position of the discontinuity ($x_d$) are also given. In all cases $\gamma = 5/3$.} 
\begin{center} 
 \begin{tabular}{rccccccccc}
 \hline
 Case & $\rho$ & $u$ & $v$ & $w$ & $p$ & $B_x$ & $B_y$ & $B_z$ & $t_{\textnormal{end}}$, $x_d$      \\ 
 \hline
 RP1 L: &  1.0    &  0.0     & 0.0    & 0.0      &  1.0     & $\frac{3}{4} \sqrt{4 \pi}$ &  $\sqrt{4 \pi}$  & 0.0   & 0.1    \\
     R: &  0.125  &  0.0     & 0.0    & 0.0      &  0.1     & $\frac{3}{4} \sqrt{4 \pi}$ & $-\sqrt{4 \pi}$  & 0.0   & 0.0    \\
 RP2 L: &  1.08   &  1.2     & 0.01   & 0.5      &  0.95    & 2.0 &  3.6     & 2.0    & 0.2         \\
     R: &  0.9891 &  -0.0131 & 0.0269 & 0.010037 &  0.97159 & 2.0 &  4.0244  & 2.0026 & -0.1        \\
 RP3 L: &  0.15   & 21.55    & 1.0    & 1.0      &  0.28    & 0.05 &  -2.0   & -1.0   & 0.04        \\
     R: &  0.1    & -26.45   & 0.0    & 0.0      &  0.1     & 0.05 &  2.0    &  1.0   &  0.0        \\
 RP4 L: &  1.0    &  0.0     & 0.0    & 0.0      &  1.0     & $1.3 \sqrt{4 \pi}$ &  $\sqrt{4 \pi}$   & 0.0   & 0.16          \\
     R: &  0.4    &  0.0     & 0.0    & 0.0      &  0.4     & $1.3 \sqrt{4 \pi}$ &  $-\sqrt{4 \pi}$  & 0.0   & 0.0           \\
 RP5 L: &  1.0    &  36.87   & -0.115 & -0.0386  &  1.0     & 4.0 &  4.0   & 1.0   & 0.03         \\
     R: &  1.0    & -36.87   & 0.0    & 0.0      &  1.0     & 4.0 &  4.0   & 1.0   & 0.0          \\
 RP6 L: &  1.7    &  0.0     & 0.0    & 0.0      &  1.7     & 3.899398 &  3.544908  & 0.0      & 0.15        \\
     R: &  0.2    &  0.0     & 0.0    & -1.496891  &  0.2   & 3.899398 &  2.785898  & 2.192064 & -0.1        \\
 \hline
 \end{tabular}
\end{center} 
 \label{tab.ic.mhd}
\end{table}

\begin{table}[!t]
 \caption{MHD equations: Conservation error for each test case and comparison of the computational efficiency between GTS and LTS algorithm using the total 
          number of element updates. } 
\begin{center} 
 \begin{tabular}{ccccccccc}
 \hline
	    & \multicolumn{5}{c}{Conservation error} & \multicolumn{3}{c}{Number of element updates} \\ 
 Case & $\rho$ & $\rho u$ & $\rho v$ & $\rho E$ & $B_y$ & GTS & LTS & GTS/LTS        \\ 
 \hline
 RP1    &  1.188E-13  &  5.066E-13  & 5.117E-13   & 2.608E-12 & 8.825E-12 & 23060      & 49400  & 2.14  \\
 RP2    &  7.159E-13  &  1.006E-12  & 1.489E-12   & 6.823E-12 & 8.284E-12 & 36313      & 50600  & 1.39  \\
 RP3    &  1.616E-12  &  3.407E-11  & 1.618E-12   & 3.888E-10 & 2.117E-11 & 33635      & 69000  & 2.05  \\
 RP4    &  5.282E-13  &  2.857E-12  & 2.356E-12   & 5.349E-12 & 1.559E-11 & 29727      & 40200  & 1.35  \\
 RP5    &  2.770E-12  &  9.040E-11  & 3.176E-12   & 1.714E-09 & 1.136E-11 & 48383      & 93600  & 1.93  \\
 RP6    &  2.139E-13  &  7.934E-13  & 3.722E-13   & 8.717E-13 & 2.140E-12 & 36234      & 70800  & 1.95  \\
 \hline
 \end{tabular}
\end{center} 
 \label{tab.eu.mhd}   
\end{table}

\paragraph{Convergence study} 

\begin{figure}
\begin{center}
\includegraphics[width=0.65\textwidth]{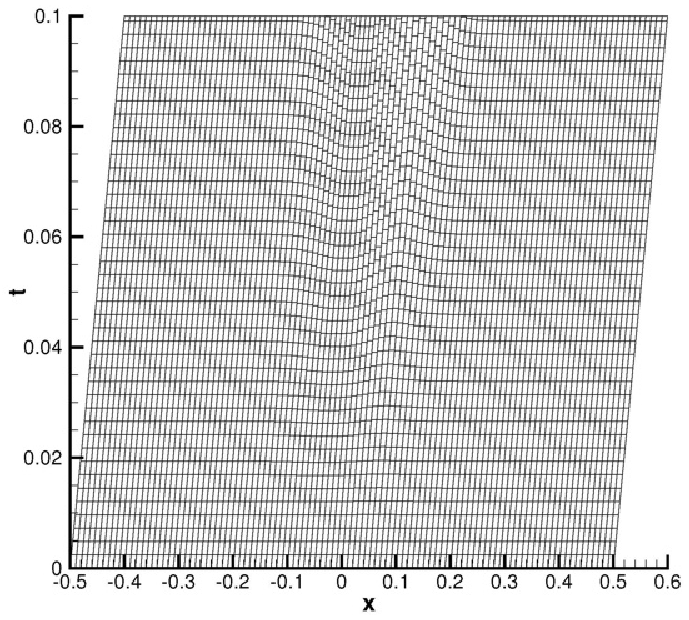}  
\caption{Space-time mesh for the numerical convergence study using a grid with 100 mesh points.}
\label{fig.convxt}
\end{center}
\end{figure}

Here we solve a smooth time-dependent problem for the MHD equations with exact solution, so that the designed order of accuracy in space and time 
of the proposed Lagrangian ADER-WENO scheme with LTS can be verified. The problem consists in a traveling Alfv\'en wave and the exact solution 
of the problem is given by $\rho(x,t)=1$, $u(x,t)=0$, $v(x,t)=1-A \exp(-\halb (x - v_a t)^2/\sigma^2)$, $w(x,t) = - \sqrt{2-v(x,t)^2}$, $p(x,t)=1$, 
$B_x(x,t)=\sqrt{4 \pi}$, $B_y(x,t)=-\sqrt{4 \pi} v(x,t)$ and $B_z(x,t)=\sqrt{4 \pi} \sqrt{2-v(x,t)^2}$, with the  Alfv\'en speed $v_a = 1$. 
The amplitude and the halfwidth of the perturbation are given by $A=0.1$ and $\sigma=0.25$, respectively. The simulation is run 
with third to fifth order schemes based on the Osher-type flux \eqref{eqn.osher}. The initial domain is $\Omega(0)=[-2;2]$. In this simulation 
the mesh velocity has been chosen as $V=v$ to get a non-trivial mesh motion. The convergence rates for variable $B_y$ at a final time of $t=0.1$
are reported in Table \ref{tab.conv}, where the number of grid cells used to discretize the domain $\Omega(t)$ is denoted by $N_G$. 
From the obtained results one can conclude that the designed high order of accuracy in space and time of the Lagrangian ADER-WENO schemes is 
maintained when the local time stepping feature (LTS) presented in this article is used. An example of a resulting space-time mesh is 
depicted for the grid $N_G=100$ in Figure \ref{fig.convxt}. 

\begin{table}  
\caption{Numerical convergence results for the ideal MHD equations using third to fifth order Lagrangian ADER-WENO finite volume schemes 
with time accurate local time stepping (LTS). The error norms refer to the variable $B_y$ (density) at time $t=0.1$.} 
\begin{center} 
\renewcommand{\arraystretch}{1.0}
\begin{tabular}{ccccccccc} 
\hline
  $N_G$ & $\epsilon_{L_2}$ & $\mathcal{O}(L_2)$ & $N_G$ & $\epsilon_{L_2}$ & $\mathcal{O}(L_2)$ & $N_G$ & $\epsilon_{L_2}$ & $\mathcal{O}(L_2)$ \\ 
\hline
  \multicolumn{3}{c}{$\mathcal{O}3$} & \multicolumn{3}{c}{$\mathcal{O}4$}  & \multicolumn{3}{c}{$\mathcal{O}5$} \\
\hline
100  & 2.1272E-04 &      & 25   & 3.1389E-03 &      & 25   & 1.5085E-03 &       \\ 
200  & 9.3217E-06 & 4.51 & 100  & 4.9217E-06 & 4.66 & 100  & 2.6843E-06 & 4.57  \\ 
400  & 1.1615E-06 & 3.00 & 200  & 2.0642E-07 & 4.58 & 150  & 3.6715E-07 & 4.91  \\ 
800  & 1.4502E-07 & 3.00 & 300  & 3.2271E-08 & 4.58 & 200  & 8.9272E-08 & 4.92  \\ 
\hline 
\hline 
\end{tabular} 
\end{center}
\label{tab.conv}
\end{table} 

\section{Conclusions} 
\label{sec.concl} 
 
To the knowledge of the author, this is the first time that a high order Lagrangian finite volume scheme with time-accurate local time stepping 
(LTS) has been presented. The design principle of the method is the use of a non-conforming space-time mesh. High order of accuracy in time
is achieved via a local space-time discontinuous Galerkin predictor, which solves element-local Cauchy problems using a weak form of the PDE in 
space-time. The initial condition of these local Cauchy problems is given by a high order WENO reconstruction. Since the cell averages 
and the spatial grid points are usually defined at different time levels, the reconstruction operator is applied on a virtual geometry and a virtual
set of cell averages that can be easily computed via $L_2$ projection using the space-time predictor. The fluxes are computed in a consistent and
conservative way using memory variables, which allow to handle the non-conforming nodes in time very easily from an algorithmic point of view. 

The algorithm has been applied to the Euler equations of compressible gas dynamics and to the ideal MHD equations. A set of 1D Riemann problems 
has been solved for both systems of conservation laws. Furthermore, numerical convergence results on a smooth unsteady test problem with exact 
solution have been shown for the MHD system using third to fifth order schemes in space and time. 

Future research will concern the extension of the local time-stepping algorithm presented in this paper to the multi-dimensional ALE ADER-WENO 
schemes presented in \cite{Lagrange3D,Lagrange2D,LagrangeNC}.

\section*{Acknowledgements}
This research has been financed by the European Research Council (ERC) under the European Union's Seventh Framework Programme 
(FP7/2007-2013) within the research project \textit{STiMulUs}, ERC Grant 
agreement no. 278267. The author would like to thank Dr. Mikhail Shashkov from LANL for the very inspiring discussions on the topic during the Oberwolfach workshop 1338b  
at the MFO Oberwolfach, Germany. 
\bibliographystyle{plain}
\bibliography{ALELTS1}

\end{document}